%% file: spdhss_h2.tex
\documentclass[11pt]{article}

\usepackage{amssymb}
\usepackage{amsmath} 
\usepackage{bbm,alltt}
\usepackage{booktabs} 
\usepackage{tabularx}
\usepackage{float} 
\usepackage{graphicx,epstopdf} 
\usepackage[caption=false]{subfig} 
\usepackage{blkarray}
\usepackage{mathtools}
\usepackage{multirow}

\usepackage{algorithm}
\usepackage{algorithmic} 

\newcommand{\algstep}[1]{\item[]\medskip\hrule\kern 2pt\hbox to \textwidth{\hspace{\labelsep}{#1}\hfill}\smallskip\hrule}

\usepackage[hang,flushmargin]{footmisc}
\makeatletter
\newcommand{\algorithmfootnote}[2][\footnotesize]{%
	\let\old@algocf@finish\@algocf@finish
	\def\@algocf@finish{\old@algocf@finish
		\leavevmode\rlap{\begin{minipage}{\linewidth}
				#1#2
		\end{minipage}}%
	}%
}
\makeatother

\makeatletter
\newcommand{\leqnomode}{\tagsleft@true\let\veqno\@@leqno}
\newcommand{\reqnomode}{\tagsleft@false\let\veqno\@@eqno}
\makeatother

\numberwithin{theorem}{section}
\numberwithin{algorithm}{section}
\allowdisplaybreaks

\usepackage{cleveref}
\usepackage{lipsum}
\usepackage{amsfonts}
\usepackage{graphicx}
\usepackage{epstopdf}
\ifpdf%
  \DeclareGraphicsExtensions{.eps,.pdf,.png,.jpg}
\else
  \DeclareGraphicsExtensions{.eps}
\fi
\usepackage{amsopn}

\title{Efficient construction of an HSS preconditioner for symmetric positive definite $\mathcal{H}^2$ matrices%
\thanks{Version of \today.}}

\author{
Xin Xing%
\thanks{Department of Mathematics, University of California, Berkeley, CA (xxing@berkeley.edu)}
\and
Hua Huang%
\thanks{School of Computational Science and Engineering, Georgia Institute of Technology,
	Atlanta, GA (huangh223@gatech.edu, echow@cc.gatech.edu).}
\and
Edmond Chow%
\footnotemark[3]
}

\begin{document}
\maketitle
\begin{abstract}
\input{doc/abstract}
\end{abstract}



\input{doc/introduction}
\input{doc/background}
\input{doc/spdhss_theorem}    
\input{doc/spdhss_construct}  
\input{doc/h2spdhss_construct}
\input{doc/numerical_test}
\input{doc/conclusion}

\bibliographystyle{plain}
\bibliography{reference}

\end{document}

%% file: doc/abstract.tex
In an iterative approach for solving linear systems with dense, ill-conditioned,
symmetric positive definite (SPD) kernel matrices, both fast
matrix-vector products and fast preconditioning operations are required.
Fast (linear-scaling) matrix-vector products are available by expressing
the kernel matrix in an $\mathcal{H}^2$ representation or an equivalent fast
multipole method representation.  
This paper is concerned with preconditioning such matrices using
the hierarchically semiseparable (HSS) matrix representation.
Previously, an algorithm was presented to construct an HSS approximation
to an SPD kernel matrix that is guaranteed to be SPD.
However, this algorithm has quadratic cost and was only designed
for recursive binary partitionings of the points defining the kernel matrix.
This paper presents a general algorithm for constructing 
an SPD HSS approximation.
Importantly, the algorithm uses the $\mathcal{H}^2$ representation
of the SPD matrix to reduce its computational complexity from quadratic
to quasilinear.
Numerical experiments illustrate how this SPD HSS approximation
performs as a preconditioner for solving linear systems
arising from a range of kernel functions.

%% file: doc/introduction.tex
\section{Introduction}\label{sec:intro}

Fast direct linear solvers exploit the hierarchical low-rank structure of matrix blocks.
This structure can be exploited in different ways (e.g., 
hierarchical off-diagonal low-rank (HODLR) \cite{ambikasaran_mathcalo_2013}, 
hierarchical semiseparable (HSS) \cite{chandrasekaran_fast_2006, xia_fast_2010}, 
recursive skeletonization \cite{ho_fast_2012}, 
hierarchical interpolative factorization (HIF) \cite{ho2016hierarchical}, 
inverse fast multipole method (IFMM) \cite{ambikasaran2014inverse, coulier2017inverse})
but, invariably, constructing these hierarchical low-rank representations
is expensive, its cost being dominated by computing accurate low-rank approximations
of matrix blocks and the associated factorizations based on these approximations.
Usually, this construction step scales superlinearly and is far more expensive than the subsequent 
solve step (which may include factorization, e.g., ULV decomposition \cite{xia_fast_2010}
for the HSS representation).

An alternative to fast direct solvers is to use iterative solvers and fast matrix-vector multiplication
provided by more general representations of the hierarchical low-rank structure (e.g.,
$\mathcal{H}$ \cite{hackbusch_sparse_1999, hackbusch_sparse_2000}, 
$\mathcal{H}^2$ \cite{hackbusch_data-sparse_2002, hackbusch_h2matrices_2000}, 
the fast multipole method (FMM) \cite{greengard_fast_1987, greengard_new_1997}, 
butterfly factorization \cite{li_butterfly_2015}).
These methods only require relatively cheap or even trivial precomputation to construct the
hierarchical low-rank representation, and can scale linearly or quasilinearly overall.
The main challenge here is slow convergence of the iterative solve for ill-conditioned matrices.

The two approaches above can be combined by using fast direct solvers as preconditioners
for the iterative solvers and using fast matrix-vector multiplication.
Referred to as \textit{rank-structured preconditioners}, the construction cost of the solvers
is greatly reduced due to the lower accuracy required of the low-rank approximations.
For symmetric positive definite (SPD) matrices, which are addressed in this paper,
it is important that the preconditioner is also SPD.  Unfortunately,
most rank-structured preconditioners, if only focusing on matrix block approximation, are not able to guarantee that
positive definiteness is preserved.

Recently, a \textit{scaling-and-compression} technique has been developed for both dense and sparse SPD matrices to
compress matrix blocks into low-rank form as part of the construction of
certain rank-structured preconditioners
\cite{cambier2020algebraic, feliu2020recursively, xia2020robust, xia_sif_2017, xin2020effectiveness, xing_spdhss_2018}.
The resulting preconditioners can be much more effective than if this
technique is not used.  It has also been found experimentally 
that preconditioners computed using this technique are more likely 
be positive definite.  
In some cases above, positive definiteness can further be guaranteed
when scaling-and-compression is used with the appropriate construction
algorithm, but the cost of constructing these SPD preconditioners is at
least quadratic for dense SPD matrices.

In this paper, we propose a quasilinear algorithm to efficiently construct an 
SPD preconditioner in HSS form by accelerating the scaling-and-compression technique,
given an $\mathcal{H}^2$ representation of the dense SPD matrix.

\begin{figure}[htp]
	\centering
	\includegraphics[width=0.7\textwidth]{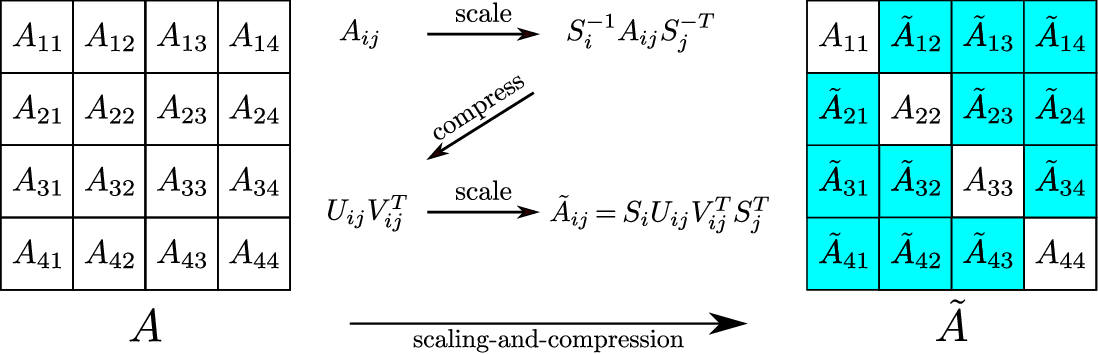}
	\caption{
		Illustration of the scaling-and-compression technique for compressing off-diagonal blocks at one level in HSS construction. 
	}
	\label{fig:intro_sc}
\end{figure}

The scaling-and-compression technique is illustrated in
\cref{fig:intro_sc} for compressing off-diagonal blocks at one level in the construction
process of an HSS representation.  The matrix $A$ is partitioned into blocks
and the compressed matrix $\tilde{A}$ is produced.  In the scaling-and-compression
technique, instead of directly compressing each off-diagonal block, the
block is scaled before compression.  Each off-diagonal block
$A_{ij}$ is scaled as $S_i^{-1} A_{ij} S_j^{-T}$, where
$S_i$ and $S_j$ are from an easily invertible symmetric factorization (e.g., Cholesky factorization) 
of the diagonal blocks,
$A_{ii} = S_i S_i^T$ and $A_{jj} = S_j S_j^T$.
The scaled off-diagonal blocks are then compressed 
into  low rank form, $S_i^{-1}A_{ij}S_j^{-T}\approx U_{ij}V_{ij}^T$.
The final low-rank approximation is $A_{ij}\approx
S_i U_{ij}V_{ij}^T S_j^T$.

The scaling step, $A_{ij} \rightarrow S_i^{-1}A_{ij}S_j^{-T}$,
requires accessing all the matrix entries and leads to quadratic
computation cost.  
If using a fixed approximation rank, the compression of all the 
scaled off-diagonal blocks, i.e., the step $S_i^{-1}A_{ij}S_j^{-T} 
\rightarrow U_{ij}V_{ij}^T$, using general algebraic methods 
such as QR decomposition and SVD also takes quadratic computation cost. 
Thus, both the scaling and compression operations in the scaling-and-compression 
technique could lead to unfavorable, quadratic HSS construction cost.

A key observation that we utilize in this paper is that if a block is
already in low-rank form, its scaling and compression can be efficiently
computed (this is utilized in \Cref{sec:bij}).  
For example, for $A_{ij} = UV^T$ with tall factors $U$ and $V$,
it is sufficient to compute and compress the two products $S_i^{-1}U$
and $S_j^{-1}V$.  If a matrix is expressed in the $\mathcal{H}^2$
representation, then the vast majority of its off-diagonal blocks is
already expressed in low-rank form.  This reduces the cost of constructing
an HSS representation that uses the scaling-and-compression technique
in the construction process.  The HSS representation generally requires
more of its off-diagonal blocks to be compressed into low-rank form than
the $\mathcal{H}^2$ representation.  However, the number of additional blocks that
need to be compressed in the HSS representation is relatively very small.
If a fixed approximation rank is used, these considerations lead 
to the quasilinear cost of the SPD HSS construction algorithm proposed in this paper.
More precisely, if a fixed rank $r$ is used for all HSS block
approximations, the new algorithm has computation dominated by $O(r\log N)$ 
matrix-vector multiplications (using the $\mathcal{H}^2$ representation)
and thus scales as $O(rN\log{N})$, where $N$ is the number of matrix
rows.  The scaling and compression of all blocks at one level of the
new construction algorithm can be performed in parallel.

There exist related ideas in the current literature.
In particular, the construction of the HSS representation and
of the butterfly factorization for a matrix can be accelerated
if a fast matrix-vector product operation is available for the matrix
\cite{ghysels2016efficient,liu2020butterfly,martinsson2011fast,rouet2016distributed}.
One application of these methods is to construct such representations for products of matrices,
where each matrix is expressed in an $\mathcal{H}$ or $\mathcal{H}^2$
representation, for example.  In a similar spirit, simple rank-structured
representations can be post-processed to construct more complicated ones,
e.g., converting an $\mathcal{H}$ representation into an $\mathcal{H}^2$
representation \cite{borm2003hierarchical}, by exploiting the efficiencies
already afforded by the existing $\mathcal{H}$ representation.

\textbf{Outline.}
Previously, a quadratic-scaling algorithm for constructing an SPD HSS 
approximation was presented~\cite{xing_spdhss_2018}.
The main concepts behind this algorithm are reviewed in \Cref{sec:spdhss}, as the new algorithm
of this paper uses the same ideas.
The earlier algorithm, however, can only construct HSS representations 
by recursively partitioning the set of matrix rows (or columns) in binary
fashion, leading to a binary partition tree (see Background, \Cref{sec:background}).
In \Cref{sec:general_construction}, we generalize the earlier algorithm to handle nonbinary partition trees.
This is a necessary step for our new algorithm because the SPD HSS
representation will be derived from an $\mathcal{H}^2$ representation
using the same partition tree, and the latter representation
can use a nonbinary partition tree.
We note that this ``generalized'' SPD HSS construction algorithm still scales quadratically.
In \Cref{sec:h2}, we propose the new algorithm that
uses an $\mathcal{H}^2$ representation of an SPD matrix
to accelerate the construction of its SPD HSS approximation,
resulting in a quasilinear algorithm.
This is the main contribution of this paper.
To demonstrate the computational cost of the new algorithm and the utility
of the SPD HSS approximation as a preconditioner,
the results of numerical experiments are shown in \Cref{sec:numerical}.

%% file: doc/background.tex
\section{Background}
\label{sec:background}

For an $N\times N$ symmetric matrix $A$, we denote its row (or column) index set
as $I = \{1,2,\ldots, N\}$.
In an applied problem, each index is associated with some element of interest,
e.g., a quadrature point, a feature vector, etc.
With a recursive partitioning of these elements of interest, the index set $I$ is partitioned 
into hierarchically enclosed subsets $\{I_i\}_{i\in\mathcal{T}}$, where 
$\mathcal{T}$ is a \textit{partition tree} that characterizes the recursive partitioning.
For each node $i\in \mathcal{T}$, $I_i$ is a subset of $I$.  
If $i$ has children $i_1, i_2, \ldots, i_m$, then $I_i = I_{i_1}\cup \cdots \cup I_{i_m}$ 
and $I_{i_a} \cap I_{i_b} = \emptyset$ for $a\neq b$. 
Often, $\mathcal{T}$ is chosen to be a binary tree, a quadtree, or an octree associated with the spatial partitioning of the elements of interest in 1-, 2-, or 3-dimensional space, respectively.
For simplicity, we assume $\mathcal{T}$ to be a perfect (fully populated in each level) $m$-ary tree.
This assumption can be lifted with minor modifications. 

The following notation is used in this paper:
\begin{itemize}
	\item For $i, j\in \mathcal{T}$, $A_{ij}$ denotes the subblock of $A$ with rows indexed by $I_i$ and columns indexed by $I_j$. 
	\item The root level of $\mathcal{T}$ is called level $L$ and the leaf level is called level 1. 
              The levels of the partition tree will be associated with levels in the hierarchical structure of a matrix.
	\item $\text{lvl}(k)$ denotes the set of nodes in level $k$ of $\mathcal{T}$. 
        \item For node $i$ in level $k$, we define $i^c = \text{lvl}(k)\setminus \{i\}$, and thus $A_{ii^c}$ denotes the off-diagonal block row of $A$ consisting of all $A_{ij}$ with $j\in i^c$.
	\item For each nonleaf node $i$, its children are denoted by $i_1, i_2, \ldots, i_m$. 
\end{itemize}

\subsection*{Low-rank approximation by projection} Given a matrix or matrix block $H\in\mathbb{R}^{n\times s}$, a general approach for compressing $H$ into rank-$r$ form is to compute a tall matrix $V\in\mathbb{R}^{n\times r}$ with orthonormal columns
whose column space, $\text{col}(V)$, is close to the \textit{principal column space} of $H$, i.e., the space spanned by the first $r$ left singular vectors of $H$. 
A rank-$r$ approximation can then be written as $H \approx V V^TH$ where $V V^T$ projects each column of $H$ onto $\text{col}(V)$. 
Such a basis matrix $V$ can be computed by SVD, QR decomposition, randomized methods, etc.

\subsection*{HSS representation}
At each level $k$, an HSS construction algorithm for a matrix $A$
compresses all the off-diagonal blocks $A_{ij}$ with $i\neq j \in \text{lvl}(k)$ into the low-rank form
\begin{equation}\label{eqn:uniform_basis}
A_{ij}\approx U_i B_{ij} U_j^T,
\end{equation}
where \textit{basis matrix} $U_i$ is shared by all the off-diagonal blocks with rows indexed by $I_i$, 
i.e., all blocks in $A_{ii^c}$, and where $U_j^T$ is similarly shared from the symmetry of $A$.
Assuming $U_i$ has orthonormal columns,
\textit{coefficient matrix} $B_{ij}$ can be computed as $U_i^T A_{ij} U_j$. Then the approximation \cref{eqn:uniform_basis} projects the columns and rows of $A_{ij}$ onto the column spaces $\text{col}(U_i)$ and $\text{col}(U_j)$, respectively.
Matrix $U_i$ captures the principal column space of $A_{ii^c}$ (to compress $A_{ii^c}$) in a recursive way.
If $i$ has children $i_1, \ldots, i_m$, then $U_i$ has the nested form
\begin{equation}\label{eqn:nested}
U_i = 
\begin{bsmallmatrix}
U_{i_1} & & \\
& \ddots & \\
& & U_{i_m}
\end{bsmallmatrix}
R_i
\end{equation}
with \textit{transfer matrix} $R_i$.
An HSS representation consists of (1) dense diagonal blocks $A_{ii}$ associated with leaf nodes and (2) low-rank representations \cref{eqn:uniform_basis} of off-diagonal blocks $A_{ij}$ at various levels that are not contained in larger off-diagonal blocks. 
Such a block $A_{ij}$ is associated with a pair of {\em sibling} nodes $i$ and $j$,
i.e., nodes $i$ and $j$ have the same parent. 
\Cref{fig:hss_mat} shows an HSS representation for a binary partition tree.

\subsection*{Recursive HSS construction}
Constructing an HSS representation starts from the leaf level (level 1) to the level below the root (level $L-1$) of $\mathcal{T}$. 
At level $1$, the original matrix $A^{(0)} = A$ has all its off-diagonal blocks $A_{ij}^{(0)}$ with $i\neq j \in \text{lvl}(1)$ compressed into the low-rank form \cref{eqn:uniform_basis} and all its diagonal blocks $A_{ii}^{(0)}$ untouched. 
This overall approximation to $A^{(0)}$ is denoted as $A^{(1)}$.
Recursively, at each level $k$, $A^{(k-1)}$ from level $(k-1)$ has its off-diagonal blocks $A_{ij}^{(k-1)}$ with $i\neq j\in\text{lvl}(k)$ compressed and is overall approximated by $A^{(k)}$. 
Lastly, $A^{(L-1)}$ is the HSS representation of $A$.

Each $A^{(k)}_{ij}$ with $i\neq j\in\text{lvl}(k)$ gives a low-rank approximation of $A_{ij}$ but is constructed \textit{indirectly} by approximating $A^{(k-1)}_{ij}$ and not the original $A_{ij}$, i.e., 
\[
A_{ij} \approx A^{(k-1)}_{ij} \approx A^{(k)}_{ij} = U_iU_i^T A^{(k-1)}_{ij} U_jU_j^T, \quad i\neq j\in\text{lvl}(k).
\]
Similarly, basis matrix $U_i$ with $i\in \text{lvl}(k)$ is constructed indirectly by compressing $A^{(k-1)}_{ii^c}$ instead of $A_{ii^c}$.
This helps enforce the nested form \cref{eqn:nested} of $U_i$.

\begin{figure}
\centering
\includegraphics[width=\textwidth]{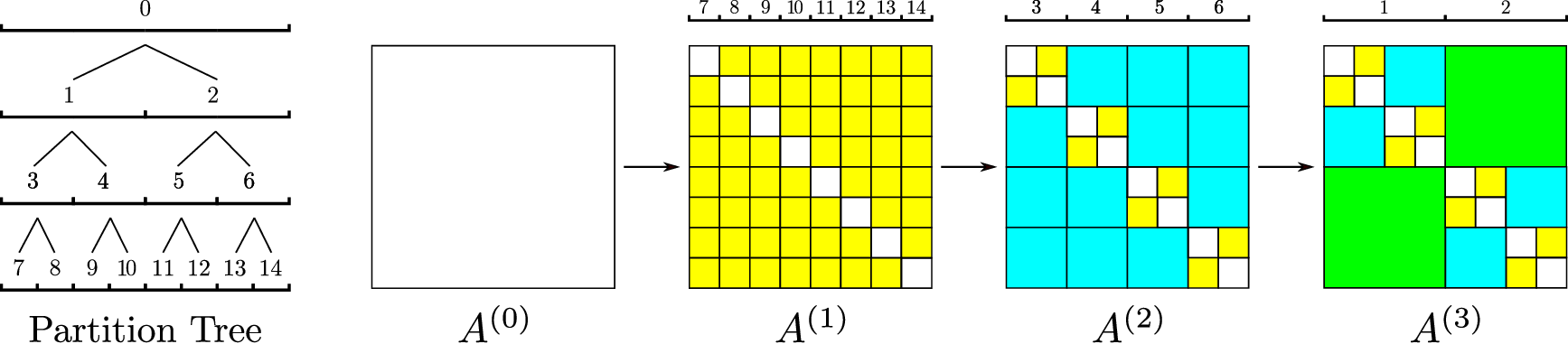}
\caption{
Recursive construction of an HSS approximation with a binary partition tree with $L=4$ levels.
The colored blocks at different levels are compressed into low-rank form, and $A^{(3)}$ gives an HSS approximation of the original matrix $A^{(0)}$.
}
\label{fig:hss_mat}
\end{figure}

The recursive HSS construction can be summarized as follows. 
For levels $k$ from $1$ to $L-1$,
\begin{multline}
   A^{(k)} = \text{diag} (\{ A_{i i}^{(k-1)} \}_{i \in \text{lvl} (k)}) \\
           + \text{diag} (\{ U_i U_i^T \}_{i \in \text{lvl} (k)}) \left[ A^{(k-1)} -
  \text{diag} (\{ A_{i i}^{(k-1)} \}_{i \in \text{lvl} (k)}) \right] \text{diag} (\{ U_i U_i^T \}_{i \in \text{lvl} (k)}) \label{eqn:recursion}
\end{multline}
where the notation $\text{diag} (\{ H_i \}_{i \in \text{lvl} (k)})$ denotes
a block diagonal matrix consisting of all blocks in $\{ H_i \}_{i \in \text{lvl} (k)}$. 
This notation will be simplified as $\text{diag} (H_i)$ with $i \in
\text{lvl}(k)$ implied by the context.
This recursive construction process is illustrated in \Cref{fig:hss_mat}.

%% file: doc/spdhss_theorem.tex
\section{Review of SPD HSS construction concepts}\label{sec:spdhss}
In this section, we review the results from Ref.~\cite{xing_spdhss_2018} that provide the cornerstone for this paper. 
Specifically, we first show how scaling-and-compression is used with the recursive HSS construction procedure
to compress the off-diagonal blocks of $A^{(k-1)}$ to obtain $A^{(k)}$ for each level $k$.
We then explain how this algorithm guarantees
that the constructed HSS approximation $A^{(L-1)}$ of $A$ is SPD.

\subsection{Scaling-and-compression technique} \label{sec:scaling-and-compression}
\begin{figure}[h]
\centering
\includegraphics[width=\textwidth]{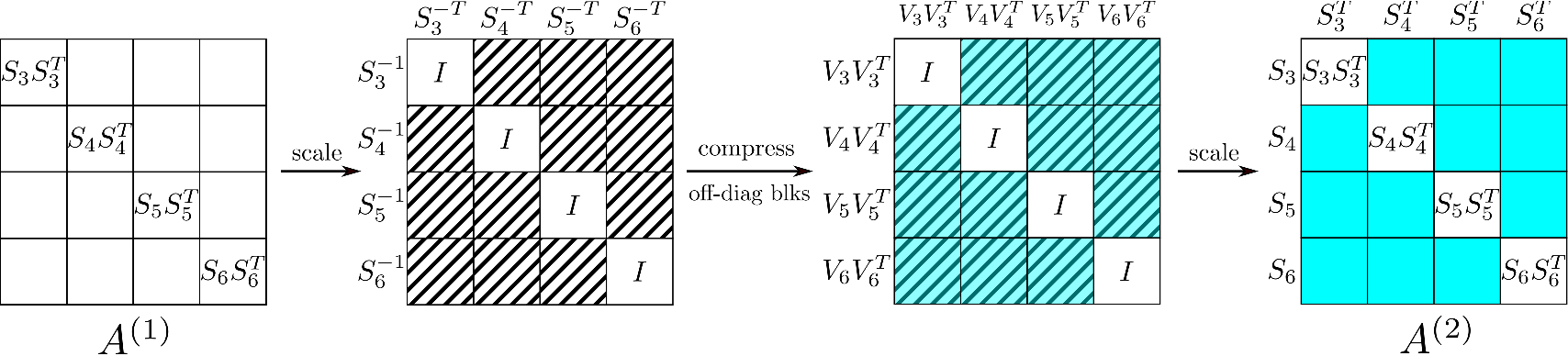}
\caption{
Illustration of the scaling-and-compression technique to compress off-diagonal blocks of $A^{(1)}$ at level $2$  to obtain $A^{(2)}$ for the example in \Cref{fig:hss_mat}. 
Note that the scaling operations are applied to all blocks, and the compression operations are only applied to off-diagonal blocks. 
}
\label{fig:scale_compress}
\end{figure}

Consider the HSS construction at level $k$ that approximates $A^{(k-1)}$ by $A^{(k)}$.
Using the scaling-and-compression technique, first compute a symmetric factorization (e.g., Cholesky decomposition) of each diagonal block $A^{(k-1)}_{ii}$ with $i\in\text{lvl}(k)$ as $A_{ii}^{(k-1)} = S_i S_i^T$.
Each off-diagonal block $A_{ij}^{(k-1)}$ with $i\neq j\in\text{lvl}(k)$ is then scaled by $S_i^{-1}$ and $S_j^{-T}$ from its left and right, respectively, as
\begin{equation}
A_{ij}^{(k-1)} \xrightarrow{\text{scale}} C_{ij}^{(k-1)} = S_i^{-1} A_{ij}^{(k-1)} S_j^{-T}.
\end{equation}
This is equivalent to multiplying $A^{(k-1)}$ by $\text{diag}(S_i^{-1})$ and $\text{diag}(S_i^{-T})$  from left and right respectively, making the diagonal blocks of $A^{(k-1)}$ be identity.
Next, compress all these scaled off-diagonal blocks $C_{ij}^{(k-1)}$.
In this paper, we use the projection approach for compression; see \cref{eqn:recursion}.  Other approaches
are possible, but may not be able to guarantee that the scaling-and-compression technique
helps give SPD approximations.
In the projection approach, compute a tall matrix $V_i$ with orthonormal columns to approximate $C_{ii^c}^{(k-1)}$ by $V_iV_i^TC_{ii^c}^{(k-1)}$ and thus compress each $C_{ij}^{(k-1)}$ as
\[
  C_{ij}^{(k-1)} \xrightarrow{\text{compress}} V_iV_i^T C_{ij}^{(k-1)} V_jV_j^T.
\]
Lastly, scale these compressed blocks back using $S_i$ and $S_j^T$ to obtain the final low-rank approximation $A_{ij}^{(k)}$ to $A_{ij}^{(k-1)}$
as 
\begin{align}
A_{ij}^{(k-1)} \approx A_{ij}^{(k)} & = S_i (V_iV_i^T C_{ij}^{(k-1)} V_jV_j^T) S_j^T \nonumber\\
 & = S_i V_i (V_i^T S_i^{-1} A_{ij}^{(k-1)} S_j^{-T} V_j) V_j^T S_j^T, \label{eqn:no_w}
\end{align}
which we write as
\begin{equation}
 A_{ij}^{(k)} = U_i B_{ij} U_j^T, \label{eqn:ubu}
\end{equation}
where we have defined the basis matrix $U_i = S_i V_i$ and the coefficient matrix
$B_{ij} = V_i^T  C_{ij}^{(k-1)} V_j
= V_i^T S_i^{-1} A_{ij}^{(k-1)} S_j^{-T} V_j$.
Thus the same notation as before is used for the basis matrix and the coefficient matrix,
regardless of whether the scaling-and-compression technique is used.
With scaling-and-compression, we again require $U_i$ to satisfy the nested form \cref{eqn:nested}.
For clarity, it is worth comparing the definition \cref{eqn:ubu} with the approximation \cref{eqn:uniform_basis}.
\Cref{fig:scale_compress} illustrates the application of the scaling-and-compression technique for compressing $A^{(1)}$ to obtain $A^{(2)}$ at level $2$ for the example of \Cref{fig:hss_mat}.

\subsection{Positive definiteness of $A^{(k)}$}
Given an SPD matrix $A$, to show that the HSS approximation $A^{(L-1)}$ constructed above is SPD, it is sufficient to show that if $A^{(k-1)}$ is SPD, then $A^{(k)}$ is also SPD. 
To begin, the low-rank approximation \cref{eqn:no_w} to an off-diagonal block can be written as
\[
A_{ij}^{(k)} = U_i W_i A_{ij}^{(k-1)} W_i^T U_i^T
\]
where we have defined $W_i = V_i^T S_i^{-1}$.
Then, the overall approximation at level $k$ is
\begin{align*}
A^{(k)} 
&= \text{diag} (A_{i i}^{(k-1)}) + \text{diag} (U_iW_i) [ A^{(k-1)} -
   \text{diag} (A_{i i}^{(k-1)}) ] \text{diag} (U_iW_i)^T \\
&= \text{diag} (U_iW_i) A^{(k-1)}\text{diag} (U_iW_i)^T + \text{diag}( A_{i i}^{(k-1)} - U_iW_iA_{i i}^{(k-1)}W_i^TU_i^T) \\
&= \text{diag} (U_iW_i) A^{(k-1)}\text{diag} (U_iW_i)^T + \text{diag}(S_i(I - V_iV_i^T)S_i^T)
\end{align*}
which shows that $A^{(k)}$ is at least positive semi-definite.
To show that $A^{(k)}$ is SPD, given that $A^{(k-1)}$ is SPD, we prove that $v^TA^{(k)}v > 0$ for any nonzero vector $v$. 
Assume  $v^TA^{(k)}v = 0$. 
Since $A^{(k-1)}$ is SPD, we have 
\[
\text{diag}(U_iW_i)^T v = 0 \quad \text{and} \quad v^T\text{diag}(S_i(I - V_iV_i^T)S_i^T)v = 0. 
\]
Let $v_i$ denote the subvector of $v$ indexed by $I_i$. The above two equations can be further written as $S_i^{-T}V_iV_i^TS_i^Tv_i = 0$ and $v_i^TS_i(I-V_iV_i^T)S_i^Tv_i = 0$, for each node $i\in\text{lvl}(k)$. 
Plugging $V_iV_i^TS_i^Tv_i = 0$ into the latter equation gives $v_i^TS_iS_i^Tv_i = 0$ which suggests $v_i = 0$. 
Thus, $A^{(k)}$ is SPD. 

%% file: doc/spdhss_construct.tex
\section{Generalized SPD HSS construction}
\label{sec:general_construction}

In the following discussion, we assume a fixed rank $r$ for the low-rank approximation of all the off-diagonal blocks in HSS construction. 
The formal construction in the previous section involves computations with large matrix blocks and leads to $O(N^3)$ computation cost. 
Its implementation with reduced, $O(N^2r)$ computation proposed in Ref.~\cite{xing_spdhss_2018} cannot be applied to nonbinary partition trees.
In this section, we generalize this quadratic SPD HSS construction method to general partition trees
and retain $O(N^2r)$ complexity.
\Cref{sec:h2} will then demonstrate how to exploit an SPD $\mathcal{H}^2$ representation to reduce the computation cost of the generalized construction method to $O(rN\log N)$.

An HSS approximation has three components: (1) diagonal blocks $A_{ii}$ for each leaf node $i$, (2) basis matrices $U_i$ for each leaf node $i$ and transfer matrices $R_i$ for each nonleaf node $i$, and (3) coefficient matrices $B_{ij}$ for each pair of siblings $i$ and $j$.
Note that although only $B_{ij}$ matrices for siblings $i$ and $j$ are used in the final HSS representation, all $B_{ij}$ matrices with any $i\neq j \in \text{lvl}(k)$ are needed during the HSS construction process. 
With fixed approximation rank $r$, the matrices $U_i$, $R_i$, and $B_{ij}$ are of dimensions $|I_i|\times r$, $mr\times r$, and $r \times r$, respectively. 

For each level $k$ from 1 to $L-1$, the following calculations are needed, for $i\neq j\in\text{lvl}(k)$:
\begin{itemize}
  \item Decomposition: $A^{(k-1)}_{i i} = S_i S_i^T$.
  
  \item Scale: $C_{ij}^{(k-1)} = S_i^{-1} A_{ij}^{(k-1)} S_j^{-T}$.
  
  \item Compute: $V_i$ to approximate $C_{ii^c}^{(k-1)}$ by $V_iV_i^T C_{ii^c}^{(k-1)}$.
  \item Compute: $B_{ij}= V_i^TC_{ij}^{(k-1)}V_j$.
  \item For leaf levels, compute: $U_i = S_i V_i$.
  \item For nonleaf levels, compute $R_i$ by solving \cref{eqn:nested},
        \begin{equation}\label{eqn:R}
        R_i = 
        \begin{bmatrix}
        V_{i_1}^T S_{i_1}^{-1} & & \\
        & \ddots & \\
        & & V_{i_m}^T S_{i_m}^{-1}
        \end{bmatrix} S_iV_i ~ .
        \end{equation}
\end{itemize}

At the leaf level ($k=1$), all the matrices in the above calculations
are small, and the components of the HSS representation, $U_i$ and $B_{ij}$,
can be computed directly with the above formulas.
At nonleaf levels, the matrices $S_i$, $C_{ij}^{(k-1)}$, and $V_i$
in the calculations can be considered large, with dimension $O(N)$ for
levels near the root of the partition tree.  However, the HSS components
actually needed at each nonleaf level are the small $r\times r$ and $mr\times r$
matrices $B_{ij}$ and $R_i$.

In the following, we show that these large matrices $S_i$,
$C_{ij}^{(k-1)}$, and $V_i$ at level $k$ can be represented using the
matrices $\{B_{ij}\}$ and $\{R_i\}$ previously computed in level $(k-1)$.
Further, $\{B_{ij}\}$ and $\{R_i\}$ at level $k$ can be computed directly using 
$\{B_{ij}\}$ and $\{R_i\}$ from level $(k-1)$.  Thus, all calculations
involving large matrix blocks can be avoided.

\paragraph{Symmetric decomposition $A^{(k-1)}_{ii} = S_i S_i^T$}
For each nonleaf node $i$ at level $k$ with children $i_1,\ldots, i_m$,
the block $A^{(k-1)}_{i i}$ can be first split as
\[
A^{(k-1)}_{ii} = 
\begin{bmatrix}
A^{(k-1)}_{i_1i_1} & \ldots & A^{(k-1)}_{i_1i_m} \\
\vdots & \ddots & \vdots \\
A^{(k-1)}_{i_mi_1} & \ldots & A^{(k-1)}_{i_mi_m}
\end{bmatrix}
= 
\begin{bmatrix}
S_{i_1}S_{i_1}^T & \ldots & U_{i_1} B_{i_1,i_m} U_{i_m}^T \\
\vdots & \ddots & \vdots \\
U_{i_m} B_{i_m,i_1} U_{i_1}^T & \ldots & S_{i_m}S_{i_m}^T
\end{bmatrix} .
\]
This matrix can then be decomposed as (using $U_{i_a} = S_{i_a}V_{i_a}$) 
\begin{equation}\label{eqn:Aii0}
A_{ii}^{(k-1)} = \begin{bmatrix}
S_{i_1} & &  \\
 & \ddots &  \\
 &  & S_{i_m}
\end{bmatrix}
\left(
I + 
\mathbf{V}_i\mathbf{B}_{ii}\mathbf{V}_i^T
\right)
\begin{bmatrix}
S_{i_1} & &  \\
 & \ddots &  \\
 &  & S_{i_m}
\end{bmatrix}^T
\end{equation}
with
\[
\mathbf{B}_{ii}= 
\begin{bmatrix}
0 & B_{i_1,i_2} & \ldots & B_{i_1,i_m} \\
B_{i_2,i_1} & 0 & \ldots & B_{i_2,i_m} \\
\vdots & \vdots & \ddots & \vdots \\
B_{i_m,i_1} & B_{i_m, i_2} & \ldots & 0
\end{bmatrix}, 
\quad 
\mathbf{V}_i = 
\begin{bmatrix}
V_{i_1} & &  \\
 & \ddots &  \\
 &  & V_{i_m}
\end{bmatrix}.
\]
We use bold typeface to denote concatenations of children blocks,
e.g., $\mathbf{B}_{ii}$ is made up of children blocks ${B_{i_a,i_b}}$ from level $(k-1)$. 

As can be verified, a symmetric factorization $I + \mathbf{V}_i\mathbf{B}_{ii}\mathbf{V}_i^T = \bar{S}_i \bar{S}_i^T$ exists with 
\begin{align}
\bar{S}_i & = I + \mathbf{V}_i((I+\mathbf{B}_{ii})^{1/2} - I)\mathbf{V}_i^T \label{eqn:Si0} \\
\bar{S}_i^{-1} & = I + \mathbf{V}_i((I+\mathbf{B}_{ii})^{-1/2} - I)\mathbf{V}_i^T \nonumber
\end{align}
which are derived from a formula in \cite{ambikasaran2014fast}.
These are the key equations that we use to generalize the SPD HSS construction method 
of Ref.~\cite{xing_spdhss_2018} for binary partition trees to nonbinary partition trees.
The positive definiteness of $A_{ii}^{(k-1)}$ guarantees the existence of $(I+\mathbf{B}_{ii})^{\pm 1/2}$. 
Matrix $\mathbf{B}_{ii}$ is of dimension $mr\times mr$ and $(I+\mathbf{B}_{ii})^{\pm 1/2}$ can be computed by the direct eigen-decomposition of $\mathbf{B}_{ii}$.
A symmetric factorization $A_{ii}^{(k-1)} = S_i S_i^T$ can be formally computed based on \cref{eqn:Aii0} and \cref{eqn:Si0} with
\begin{equation}\label{eqn:Si}
S_i = 
\begin{bmatrix}
S_{i_1} & &  \\
 & \ddots &  \\
 &  & S_{i_m}
\end{bmatrix}
\bar{S}_i.
\end{equation}

\paragraph{Scaled off-diagonal blocks $C_{ij}^{(k-1)} = S_i^{-1} A_{ij}^{(k-1)} S_j^{-T}$} 
For nonleaf nodes $i\neq j$ at level $k$ with children $i_1,\ldots, i_m$ and $j_1, \ldots, j_m$, the quantity $A_{ij}^{(k-1)}$ can be written as
\[
\begin{bmatrix}
A^{(k-1)}_{i_1j_1} & \ldots & A^{(k-1)}_{i_1j_m} \\
\vdots & \ddots & \vdots \\
A^{(k-1)}_{i_mj_1} & \ldots & A^{(k-1)}_{i_mj_m}
\end{bmatrix}
\!\!=\!\!
\begin{bmatrix}
U_{i_1} & &  \\
 & \!\! \ddots \!\! &  \\
 &  & U_{i_m}
\end{bmatrix}\!\!
\begin{bmatrix}
B_{i_1j_1} & \ldots & B_{i_1j_m} \\
\vdots & \ddots & \vdots \\
B_{i_mj_1} & \ldots & B_{i_mj_m}
\end{bmatrix}\!\!
\begin{bmatrix}
U_{j_1} & &  \\
 & \!\! \ddots \!\! &  \\
 &  & U_{j_m}
\end{bmatrix}^T
\]
where the middle matrix is denoted as $\mathbf{B}_{ij}\in\mathbb{R}^{mr\times mr}$. 
Noting its difference from $B_{ij}$, this bold typeface $\mathbf{B}_{ij}$ consists of children blocks $B_{i_a,j_b}$ in level $(k-1)$. 
Based on \cref{eqn:Si0}, \cref{eqn:Si}, $U_i = S_iV_i$, and the above equation, the scaled block $C_{ij}^{(k-1)}$ by definition can be computed as 
\begin{align}
C_{ij}^{(k-1)} 
 = S_i^{-1} A_{ij}^{(k-1)} S_j^{-T} &= \bar{S}_i^{-1}\mathbf{V}_i \mathbf{B}_{ij} \mathbf{V}_j^T \bar{S}_j^{-T} \nonumber\\
& = \mathbf{V}_i(I + \mathbf{B}_{ii})^{-1/2} \mathbf{B}_{ij} (I + \mathbf{B}_{jj})^{-1/2}\mathbf{V}_j^T. \label{eqn:Cij}
\end{align}

\paragraph{Calculation of $V_i$} 
Recall that we desire $V_i$ such that $\text{col}(V_i)$ approximates $\text{col}(C_{ii^c}^{(k-1)})$.
Further, $V_i$ must satisfy $\text{col}(V_i) \subset \text{col}(C_{ii^c}^{(k-1)})$ in order to guarantee the nested form of $U_i$ in \cref{eqn:nested}; see Ref.~\cite{xing_spdhss_2018}.
From \cref{eqn:Cij}, each scaled block $C_{ij}^{(k-1)}$ has its column space contained in $\text{col}(\mathbf{V}_i)$.
Thus, $V_i$ can be represented by
\begin{equation}\label{eqn:Vi}
V_i 
= \mathbf{V}_i \bar{V}_i 
= 
\begin{bmatrix}
V_{i_1} & &  \\
 & \ddots &  \\
 &  & V_{i_m}
\end{bmatrix}
\bar{V}_i,
\end{equation}
where the small matrix $\bar{V}_i \in \mathbb{R}^{mr\times r}$ is computed
with orthonormal columns to minimize the error of the required approximation $C_{ii^c}^{(k-1)} \approx V_i V_i^T C_{ii^c}^{(k-1)}$. 
Noting that all $V_{i_a}$ blocks have orthonormal columns and using \cref{eqn:Cij} and \cref{eqn:Vi}, the minimization problem can be converted as
\begin{equation}\label{eqn:minimization}
\min_{V_i} \|C_{ii^c}^{(k-1)} - V_i V_i^T C_{ii^c}^{(k-1)}\|_F = \min_{\bar{V_i}} \left\| E_{ii^c} - \bar{V}_i \bar{V}_i^T E_{ii^c} \right\|_F,
\end{equation}
where $E_{ii^c}$ is the horizontal concatenation of all blocks $(I + \mathbf{B}_{ii})^{-1/2} \mathbf{B}_{ij} (I + \mathbf{B}_{jj})^{-1/2}$
with nodes $j\in i^c$. 
We note that $E_{ii^c}$ is a small matrix of dimension $mr\times (|\text{lvl}(k)|-1)r$.
Thus, $\bar{V}_i$ can be directly computed to capture the principal column space of $E_{ii^c}$. 
In \Cref{sec:calc_v}, we will discuss how to more efficiently compute $\bar{V}_i$.

\paragraph{Calculation of $B_{ij}=V_i^TC_{ij}^{(k-1)}V_j$}
Based on the calculations of $C_{ij}^{(k-1)}$ in \cref{eqn:Cij} and $V_i$ in \cref{eqn:Vi}, $B_{ij}$ can be directly computed as
\begin{align}
B_{ij}
& = 
\left(
\bar{V}_i^T
\mathbf{V}_i^T
\right)
\mathbf{V}_i(I + \mathbf{B}_{ii})^{-1/2} \mathbf{B}_{ij} (I + \mathbf{B}_{jj})^{-1/2}\mathbf{V}_j^T
\left(
\mathbf{V}_j \bar{V}_j
\right) \nonumber \\
& = 
\bar{V}_i^T (I + \mathbf{B}_{ii})^{-1/2} \mathbf{B}_{ij}(I + \mathbf{B}_{jj})^{-1/2}\bar{V}_j
\label{eqn:Mijk}
\end{align}
where all matrices in the second equation are of small dimensions.
In \Cref{sec:bij}, we will discuss how to reduce the number of matrices $B_{ij}$ that
need to be computed.

\paragraph{Calculation of $R_i$}
Based on the above calculation of $S_i$ in \cref{eqn:Si} and $V_i$ in \cref{eqn:Vi}, $R_i$ defined by \cref{eqn:R} can be directly computed as
\begin{align}
R_i 
& = 
\begin{bmatrix}
V_{i_1}^T S_{i_1}^{-1} & & \\
& \ddots & \\
& & V_{i_m}^T S_{i_m}^{-1}
\end{bmatrix}
\begin{bmatrix}
S_{i_1} & & \\
& \ddots & \\
& & S_{i_m}
\end{bmatrix}
\bar{S}_i
\begin{bmatrix}
V_{i_1} & & \\
& \ddots & \\
& & V_{i_m}
\end{bmatrix}
\bar{V}_i \nonumber \\
& = 
(I+\mathbf{B}_{ii})^{1/2} \bar{V}_i \label{eqn:Rik}
\end{align}
where, again, all matrices in the second equation are of small dimensions.

To summarize, the actual computations needed at level $k$ include the calculation of
$(I+\mathbf{B}_{ii})^{\pm 1/2}$ in \cref{eqn:Si0}, $\bar{V}_i$ in \cref{eqn:minimization}, $B_{ij}$ in
\cref{eqn:Mijk}, and $R_i$ in \cref{eqn:Rik}.
The pseudocode of this generalized SPD HSS construction process based on scaling-and-compression is shown in \Cref{alg:spdhss}. 
The overall computation and peak storage costs of \Cref{alg:spdhss} are both $O(N^2r)$ and the
constructed SPD HSS representation has $O(Nr)$ storage cost.

\begin{algorithm}
\caption{Generalized SPD HSS construction}
\label{alg:spdhss}

\begin{algorithmic}
\REQUIRE HSS rank $r$, an SPD matrix $A$
\ENSURE an SPD HSS approximation with $\{A_{ii}\}, \{B_{ij}\}, \{U_i\}, \{R_i\}$

\textbf{At the leaf level}

\STATE \quad compute the Cholesky decomposition $A_{i i} = S_i S_i^T$, $\forall i \in \text{lvl}(1)$

\STATE \quad compute the scaled off-diagonal block $C_{i j}^{(0)} = S_i^{-1} 
		A_{i j} S_j^{- T}$, $\forall i\neq j \in \text{lvl}(1)$

\STATE \quad compute $V_i\in \mathbb{R}^{|I_i| \times r}$ satisfying 
\STATE \qquad -- $V_i$ has orthonormal columns and  $\text{col} (V_i) \subset \text{col} (C_{i i^c}^{(0)})$
\STATE \qquad -- $V_i$ should minimize $\| C_{i i^c}^{(0)}
- V_i V_i^T C_{i i^c}^{(0)} \|_F$
\STATE \quad compute $B_{i j}= V_i^T C_{i j}^{(0)}
V_j$, $\forall i\neq j \in \text{lvl}(1)$

\STATE \quad set $U_i = S_i V_i$, $\forall i \in \text{lvl}(1)$

\FOR {$k = 2,3,\ldots,L-1$}

\STATE compute $(I + \mathbf{B}_{ii})^{\pm 1/2}$ via the eigen-decomposition of $\mathbf{B}_{ii}$, $\forall i \in \text{lvl}(k)$

\STATE compute $(I + \mathbf{B}_{ii})^{-1/2} \mathbf{B}_{ij} (I + \mathbf{B}_{jj})^{-1/2}$, $\forall i\neq j \in \text{lvl}(k)$
\STATE assemble $E_{ii^c}$ in \cref{eqn:minimization} and compute $\bar{V}_i$ satisfying 

\STATE \quad -- $\bar{V}_i$ has orthonormal columns and $\text{col} (\bar{V}_i) \subset 
\text{col} (E_{ii^c})$

\STATE \quad -- $\bar{V}_i$ should minimize 
$\| E_{ii^c} - \bar{V}_i\bar{V}_i^T E_{ii^c} \|_F$

\STATE compute $B_{i j} = \bar{V}_i^T (I + \mathbf{B}_{ii})^{-1/2} \mathbf{B}_{ij} (I + \mathbf{B}_{jj})^{-1/2}\bar{V}_j$,  $\forall i \neq j \in \text{lvl}(k)$

\STATE set $R_i = (I+\mathbf{B}_{ii})^{1/2} \bar{V}_i$, $\forall i \in \text{lvl}(k)$
\ENDFOR
\end{algorithmic}
\end{algorithm}

%% file: doc/h2spdhss_construct.tex
\section{Accelerated SPD HSS construction with quasilinear computation}
\label{sec:h2}
The generalized SPD HSS construction algorithm of the previous section
has quadratic computation cost.  In this section, we show how to
reduce the cost to quasilinear if we can utilize an existing $\mathcal{H}^2$ representation
of the SPD matrix.  Below, we first give 
necessary background on $\mathcal{H}^2$ representations.

\subsection{$\mathcal{H}^2$ representation}
Like the HSS representation,
the $\mathcal{H}^2$ representation of a matrix $A$ is based on a partition tree $\mathcal{T}$
and a hierarchical index set $\{I_i\}_{i\in\mathcal{T}}$.
For each node $i$ at each level $k$ of a partition tree, we define
a node set $\mathcal{F}_i\subset\text{lvl}(k)$ that contains all the nodes in level $k$ that are in
the ``far field'' of node $i$.
More precisely, if the indices are associated with points in space,
$\mathcal{F}_i$ can be defined as the set of nodes $j\in\text{lvl}(k)$ such that
the points associated with $I_j$ are well separated from the points associated with $I_i$.
In particular, the HSS representation is a specific $\mathcal{H}^2$ representation with $\mathcal{F}_i = i^c = \text{lvl}(k)\setminus \{i\}$.

The node set $\mathcal{F}_i$
specifies which blocks will be compressed in the $\mathcal{H}^2$ representation of $A$.
At each level $k$, all blocks $A_{ij}$ with $i\in\text{lvl}(k)$ and $j\in \mathcal{F}_i\subset\text{lvl}(k)$ are compressed as
\begin{equation}\label{eqn:uniformbasis_h2}
A_{ij} = U_i^{\mathcal{H}^2} B_{ij}^{\mathcal{H}^2} (U_j^{\mathcal{H}^2})^T
\end{equation}
assuming that the $\mathcal{H}^2$ representation is exact.
Here we use the superscript ``$\mathcal{H}^2$'' to distinguish the corresponding components of $\mathcal{H}^2$ from those of HSS. 
Like for HSS, the basis matrix $U_i^{\mathcal{H}^2}$ is shared by all blocks $A_{ij}$ with $j\in\mathcal{F}_i$ and is computed to capture the principal column space of $A_{i\mathcal{F}_i} = [A_{ij}]_{j\in\mathcal{F}_i}$ (corresponding to $A_{ii^c}$ in HSS).
Further, $U_i^{\mathcal{H}^2}$ satisfies the nested form \cref{eqn:nested} as well.

\begin{figure}[htb]
\centering
\includegraphics[width=\textwidth]{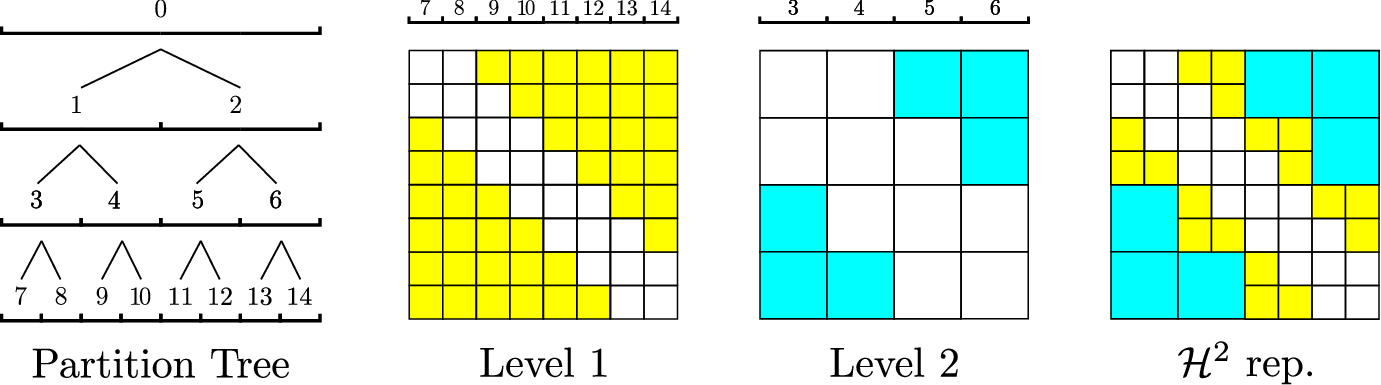}
\caption{
Illustration of an $\mathcal{H}^2$ representation with a binary partition tree.
The colored blocks $A_{ij}$ at different levels  satisfy $j\in\mathcal{F}_i$ and are compressed into low-rank form.
For each node $i$, this example defines $\mathcal{F}_i$ as the set of nodes in the same level that are not adjacent to $i$. 
}
\label{fig:h2_mat}
\end{figure}

An $\mathcal{H}^2$ representation consists of (1) dense blocks $A_{ij}$ with $j\notin\mathcal{F}_i$ at the leaf level and (2) low-rank representations \cref{eqn:uniformbasis_h2} of blocks $A_{ij}$ with $j\in\mathcal{F}_i$ at various levels that are not contained in larger low-rank blocks. 
Such a low-rank block $A_{ij}$ is associated with $i,j$ satisfying the condition $j\in\mathcal{F}_i$ but $\text{par}(j) \notin \mathcal{F}_{\text{par}(i)}$ ($\text{par}(i)$ denotes the parent of $i$).
\Cref{fig:h2_mat} gives an illustration of an $\mathcal{H}^2$ matrix with a binary partition tree. 

In practical problems, a proper definition of $\mathcal{F}_i$ can guarantee that all compressed blocks $A_{i\mathcal{F}_i}$ have numerical ranks bounded by a small constant independent of the matrix size and thus the $\mathcal{H}^2$ representation can have linear-scaling matrix-vector multiplications.
In the case of HSS, by defining $\mathcal{F}_i = i^c$, the maximum numerical rank of all $A_{i\mathcal{F}_i}$ blocks usually increases with the matrix size, and thus leads to superlinear complexities in HSS construction and other HSS computations.

\subsection{Quasilinear SPD HSS construction}

%

In the generalized SPD HSS construction algorithm of \Cref{sec:general_construction}, the computation and storage costs are dominated by those related
to the coefficient matrices, $B_{ij}$.
In each level $k$, there are $|\text{lvl}(k)|(|\text{lvl}(k)|-1)$ such matrices. 
Each $B_{ij}$ is computed recursively from the leaf level to level $k$ using \cref{eqn:Mijk} and is ultimately computed from the original matrix block $A_{ij}$. 

Before proceeding, we define $\Phi_i$, which will be used in this section.
Observe that, at level $k$,
\[
B_{ij} = V_i^T S_i^{-1} A_{ij}^{(k-1)} S_j^{-T} V_j, \quad i\neq j \in \text{lvl}(k)
\]
and that $B_{ij}$ is computed recursively by applying multiple matrices to $A_{ij}$ on its left and right.
To emphasize this relationship, we define $\Phi_i$ such that
\begin{align}
B_{ij} & = 
\Phi_i A_{ij} \Phi_j^T, \quad i\neq j \in \text{lvl}(k) 
\label{eqn:bij_hss} \\
\Phi_i & = 
\left\{
\begin{array}{ll}
 V_i^TS_i^{-1} &  i \text{ is a leaf node}\\
\bar{V}_i^T (I + \mathbf{B}_{ii})^{-1/2}\begin{bsmallmatrix}
\Phi_{i_1} & & \\
& \ddots & \\
& & \Phi_{i_m}
\end{bsmallmatrix} & i \text{ has children } i_1, \ldots, i_m
\end{array}
\right. ,
\label{eqn:phi}
\end{align}
where this nested representation of $\Phi_i$ is derived from \cref{eqn:Mijk}.

To compute the matrices $B_{ij}$, the matrices $V_i$ for leaf nodes $i$ and
$\bar{V}_i$ for nonleaf nodes $i$ are needed.  We discuss how
$V_i$ and $\bar{V}_i$ are computed using a randomized algorithm in
\Cref{sec:calc_v}.  The matrix-vector products required in the
randomized algorithm are performed efficiently using an
$\mathcal{H}^2$ representation of the SPD matrix that we assume
to be available.

To compute a $B_{ij}$ matrix at level $k$, the $B_{ij}$ matrices that
are needed from lower levels may already have a low-rank form in
the $\mathcal{H}^2$ representation.  Thus, the recursion for computing
$B_{ij}$ at level $k$ can stop and does not need to proceed to
the leaf level.  We discuss this in \Cref{sec:bij}.

\subsubsection{Calculation of $V_i$ and $\bar{V}_i$} 
\label{sec:calc_v}

In the generalized SPD HSS construction algorithm (\Cref{alg:spdhss}),
$V_i$ and $\bar{V}_i$ are computed as follows.
At leaf nodes $i$, the matrix $V_i$ is computed to approximate
$C_{i i^c}^{(0)}$ by $V_i V_i^T C_{i i^c}^{(0)}$ with
the constraint that $V_i$ has orthonormal columns and
$\text{col} (V_i) \subset \text{col} (C_{i i^c}^{(0)})$.
At nonleaf nodes $i$, the matrix $\bar{V}_i$ is computed to approximate
$E_{ii^c}$ by $\bar{V}_i\bar{V}_i^T E_{ii^c}$ with
the constraint that $\bar{V}_i$ has orthonormal columns and
$\text{col} (\bar{V}_i) \subset \text{col} (E_{ii^c})$.

For the accelerated SPD HSS algorithm, we will compute
$V_i$ and $\bar{V}_i$ using a randomized algorithm \cite{halko_finding_2011}.
(For completeness, we give the randomized algorithm in \cref{alg:matvei_random}.)
However, instead of using matrix-vector products with
$C_{i i^c}^{(0)}$ and $E_{ii^c}$, respectively, which would be the standard
approach, we will use matrix-vector products with
alternative matrices that have almost the same column spaces
as $C_{i i^c}^{(0)}$ and $E_{ii^c}$, respectively, to reduce cost.

To see what alternative matrices we propose using, we first write
the matrices $C_{i i^c}^{(0)}$ and $E_{ii^c}$ explicitly as
\begin{align*}
C_{ii^c}^{(0)} 
& = \left[S_i^{-1}A_{ij}S_j^{-T}\right]_{j\in i^c} \\
& = S_i^{-1} A_{ii^c} \, \text{diag}(\{S_j^{-T}\}_{j\in i^c}), \\
E_{ii^c} &= \left[(I+\mathbf{B}_{ii})^{-1/2} \mathbf{B}_{ij}(I+\mathbf{B}_{jj})^{-1/2}\right]_{j\in i^c}\\
& = (I+\mathbf{B}_{ii})^{-1/2} 
\left\{
\begin{bsmallmatrix}
B_{i_1 j_1} & \cdots & B{i_1, j_m} \\
\vdots           & \ddots & \vdots \\
B_{i_mj_1} & \cdots & B_{i_mj_m}
\end{bsmallmatrix}
\right\}_{j\in i^c}
\, \text{diag}(\{(I+\mathbf{B}_{jj})^{-1/2}\}_{j\in i^c})\\
&= (I+\mathbf{B}_{ii})^{-1/2} \left[\begin{smallmatrix}
\Phi_{i_1} & & \\
& \ddots & \\
& & \Phi_{i_m}
\end{smallmatrix}\right] A_{ii^c} \,
\text{diag}(\Phi_{j_s}) \, \text{diag}(\{(I+\mathbf{B}_{jj})^{-1/2}\}_{j\in i^c})
\end{align*}
where $\text{diag}(\Phi_{j_s})$ denotes the block diagonal matrix made up of all $\Phi_{j_s}$ with $j_s$
being a child of any node $j\in i^c$.
The last equation above is from substituting 
\cref{eqn:bij_hss} into $B_{i_a j_b}$ of its previous equation.

Instead of approximating the column spaces of $C_{i i^c}^{(0)}$ (when $i$ is a leaf node)
and $E_{ii^c}$ (when $i$ is a nonleaf node),
we approximate the column spaces of $\Lambda_{ii^c}$, defined as,
\[
\Lambda_{ii^c} = 
\left\{
\begin{array}{ll}
S_i^{-1}A_{ii^c} & i \text{ is a leaf node} \\
(I+\mathbf{B}_{ii})^{-1/2} \left[\begin{smallmatrix}
\Phi_{i_1} & & \\
& \ddots & \\
& & \Phi_{i_m}
\end{smallmatrix}\right]A_{ii^c} & i \text{ has children } i_1, \ldots, i_m
\end{array}
\right. .
\]
Block $\Lambda_{ii^c}$ differs from $C_{ii^c}^{(0)}$ and $E_{ii^c}$ in that there is no matrix applied to the right of $A_{ii^c}$. 
This choice of $\Lambda_{ii^c}$ is for the efficiency of computing the corresponding matrix-vector products in the randomized algorithm. 
It is theoretically possible that computing $V_i$ and $\bar{V}_i$ using $\Lambda_{ii^c}$ may affect the approximation accuracy of $C_{ii^c}^{(0)}\approx V_iV_i^TC_{ii^c}^{(0)}$ and $E_{ii^c}\approx V_iV_i^TE_{ii^c}$. 
Since our goal is to construct a low-accuracy SPD HSS preconditioner, this possible slight deterioration of the approximation accuracy may be tolerable. 

\begin{algorithm}
	\caption{Randomized algorithm for computing a basis matrix $U$ for $H$}
	\label{alg:matvei_random}
	\begin{algorithmic}
		\REQUIRE{$H \in \mathbb{R}^{n\times m}$, rank $r$, over-sampling parameter $p$}
		\ENSURE{$U$ from the rank-$r$ approximation $H \approx UU^T H$ with $U^TU=I$
		}
		\STATE \textbf{Step 1:} Generate an $m \times (p+r)$ random matrix $\Omega$ whose entries follow the standard normal distribution 
		
		\STATE \textbf{Step 2:} Compute $\Psi = H \Omega$ 
		
		\STATE \textbf{Step 3:} Compute the pivoted QR decomposition $\Psi P = QR$, and set $U$ to be the first $r$ columns of $Q$ 
	\end{algorithmic}
\end{algorithm}

The product of $\Lambda_{ii^c}$ and a block of random vectors
involves first computing the product of $A_{ii^c}$ with random vectors.
Thus we first compute the products,
\begin{equation}\label{eqn:Yk}
Y^{(k)} = (A - \text{diag}(\{A_{ii}\}_{i\in\text{lvl}(k)})) \, \Omega, \quad k = 1,2,\ldots, L-1
\end{equation}
(one for each nonroot level)
where $\Omega \in \mathbb{R}^{N\times (r+p)}$ is a random matrix,
given that we desire rank $r$ approximations using an oversampling
parameter $p$.  The quantity in the outer brackets of \cref{eqn:Yk} is just the matrix
$A$ without its block diagonal part at each level $k$.  These products $Y^{(k)}$
can be computed efficiently using the $\mathcal{H}^2$ representation of
$A$ and just neglecting the block diagonal parts during multiplication.
The desired products $A_{ii^c} \Omega_{i^c}$, where $\Omega_{i^c}$ denotes
the row subset of $\Omega$ associated with $i^c$, can be extracted as
the row subsets of $Y^{(k)}$ associated with each $i\in\text{lvl}(k)$,
denoted by $Y_i^{(k)}$.

To complete the multiplication by $\Lambda_{ii^c}$, we now apply $S_i^{-1}$ (if $i$ is a leaf node) or 
$(I+\mathbf{B}_{ii})^{-1/2} \text{diag}(\{\Phi_{i_1}, \ldots, \Phi_{i_m}\})$ (if $i$ is a nonleaf node)
to $Y_i^{(k)}$ to obtain the product $\Lambda_{ii^c} \Omega_{i^c}$
needed in step 2 of \Cref{alg:matvei_random}. 
The product $S_i^{-1}Y_i^{(1)}$ for each leaf node can be directly computed. 
The product $(I+\mathbf{B}_{ii})^{-1/2}\text{diag}(\{\Phi_{i_1}, \ldots, \Phi_{i_m}\})Y_i^{(k)}$ for each nonleaf node at level $k$ needs to be recursively computed from level 1 to level $(k-1)$, since $\Phi_{i_s}$ is recursively defined in \cref{eqn:phi}.
This recursive computation can be unfolded into local computations at each descendant of node $i$ from level 1 to level $(k-1)$ as shown in \cref{alg:spdhss_matvec}.

The complexity of computing each $Y^{(k)}$ in \cref{eqn:Yk} is
$O((r+p)N)$ due to the linear-scaling of $\mathcal{H}^2$ matrix-vector multiplication.
Since there are a logarithmic number of levels, the overall complexity
for the randomized algorithm is $O(rN\log{N})$ for both computation and storage,
assuming $p$ is a small constant.
The cost of the pivoted QR decompositions in the randomized algorithm is small
because $\Lambda_{ii^c} \Omega_{i^c}$ is a small matrix of dimension $|I_i|\times
(r+p)$ for a leaf node and $mr\times (r+p)$ for a nonleaf node.

\begin{algorithm}
\caption{Level-by-level computation of the special products}
\label{alg:spdhss_matvec}

\begin{algorithmic}
\REQUIRE $\{S_i, V_i\}$ at level 1, $\{\mathbf{B}_{ii}, \bar{V}_i\}$ at levels $2, \ldots, k-1$, $\{\mathbf{B}_{ii}\}$ at level $k$
\ENSURE $(I+\mathbf{B}_{ii})^{-1/2}\text{diag}(\{\Phi_{i_1}, \ldots, \Phi_{i_m}\})Y_i^{(k)}$ for each $i\in\text{lvl}(k)$

\textbf{At the leaf level}
	\STATE \quad Let $T_i = Y_i^{(k)}$ for each $i\in\text{lvl}(1)$ be the row subset of $Y^{(k)}$ indexed by  $I_i$
	\STATE \quad Compute $T_i = V_i^TS_i^{-1} T_i$ (in-place computation)
\FOR {$l = 2,3,\ldots,k-1$}
	\STATE Let $T_i = [T_{i_1}^T, \ldots, T_{i_m}^T]^T, \forall i\in\text{lvl}(l)$ be the vertical concatenation of $T_{i_1}, \ldots, T_{i_m}$
	\STATE Compute $T_i = \bar{V}_i^T(I+\mathbf{B}_{ii})^{-1/2} T_i$
\ENDFOR

Define $(I+\mathbf{B}_{ii})^{-1/2}[T_{i_1}^T, \ldots, T_{i_m}^T]^T$ for each $i\in\text{lvl}(k)$ as the output
\end{algorithmic}
\end{algorithm}

\subsubsection{Calculation of $B_{ij}$}
\label{sec:bij}

%

\begin{figure}[htb]
	\centering
	\includegraphics[width=\textwidth]{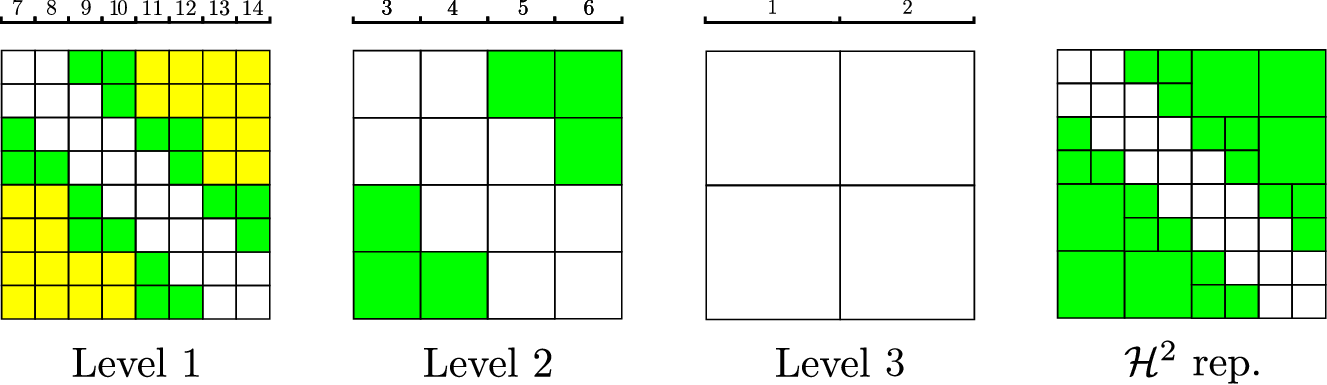}
	\caption{
		Illustration of three types of blocks at each partition level of the matrix from \Cref{fig:h2_mat}. 
		The white blocks are Type-1, the yellow blocks are Type-2, and the green blocks are Type-3. 
		An $\mathcal{H}^2$ representation is made up of Type-1 blocks from the leaf level and Type-3 blocks from all levels. 
	}
	\label{fig:h2_mat_type}
\end{figure}

We first define some nomenclature for the blocks in an $\mathcal{H}^2$ representation.  
At each level $k$, we categorize all the blocks $\{A_{ij}\}$ with
$i,j\in\text{lvl}(k)$ into three types as follows.  The colors for each
type refer to the colors in \Cref{fig:h2_mat_type} which illustrates the categorization.
\begin{itemize}
\item Type-1 (white): $j\notin\mathcal{F}_i$. 
\item Type-2 (yellow): $j\in\mathcal{F}_i$ and $A_{ij}$ is contained in a larger low-rank block at some upper level, i.e., $\text{par}(j) \in \mathcal{F}_{\text{par}(i)}$. 
\item Type-3 (green): $j\in\mathcal{F}_i$ and $A_{ij}$ is represented in low-rank form, i.e., $\text{par}(j) \notin \mathcal{F}_{\text{par}(i)}$. 
\end{itemize}
Type-1 blocks are either stored in dense form or consist of Type-1 and Type-3 blocks at next
the lower level. 
Type-2 blocks are contained in larger Type-3 blocks.

Recall from \cref{eqn:Mijk} that, in the generalized SPD HSS construction algorithm,
$B_{ij}$ at a nonleaf level is recursively calculated using
\[
B_{ij}
= \bar{V}_i^T (I + \mathbf{B}_{ii})^{-1/2} \mathbf{B}_{ij}(I + \mathbf{B}_{jj})^{-1/2}\bar{V}_j .
\]
We now discuss how the $\mathcal{H}^2$ representation can be used to
reduce the number of $B_{ij}$ matrices that need to be calculated in HSS construction.
There are three cases, corresponding to the three types of blocks.

\paragraph{Case 1}
If we need the HSS coefficient matrix $B_{ij} = \Phi_i A_{ij} \Phi_j^T$ 
and the corresponding $A_{ij}$ is a Type-3 block in the $\mathcal{H}^2$ representation, i.e.,
\[
A_{ij} = U_i^{\mathcal{H}^2} B_{ij}^{\mathcal{H}^2} (U_j^{\mathcal{H}^2})^T
\]
(from \cref{eqn:uniformbasis_h2}), then $B_{ij}$ can be computed as 
\begin{equation}\label{eqn:type3_Mij}
B_{ij}=(\Phi_i U_i^{\mathcal{H}^2}) B_{ij}^{\mathcal{H}^2} (\Phi_jU_j^{\mathcal{H}^2} )^T. 
\end{equation}
Thus, instead of recursively computing $B_{ij}$, we can first compute
$\Phi_i U_i^{\mathcal{H}^2}$ for each $i$ and use \cref{eqn:type3_Mij} to immediately compute $B_{ij}$, which only contains products of small matrices.
By this approach, we do not have to compute any $B_{i_aj_b}$ for the descendants $i_a$ of $i$ and $j_b$ of $j$ at lower levels. 
Any $B_{ij}$ in the HSS representation that corresponds to a Type-2 block $A_{ij}$ in the
$\mathcal{H}^2$ representation is no longer needed since Type-2 blocks are enclosed in Type-3 blocks.

For each node $i$, $\Phi_i U_i^{\mathcal{H}^2}$ can be recursively computed as (utilizing \cref{eqn:nested}),
\begin{align*}
\Phi_i U_i^{\mathcal{H}^2} 
& = 
\bar{V}_i^T (I + \mathbf{B}_{ii})^{-1/2}
\begin{bsmallmatrix}
\Phi_{i_1} & & \\
& \ddots & \\
& & \Phi_{i_m}
\end{bsmallmatrix}
\begin{bsmallmatrix}
U_{i_1}^{\mathcal{H}^2}  & & \\
& \ddots & \\
& & U_{i_m}^{\mathcal{H}^2} 
\end{bsmallmatrix}
R_i^{\mathcal{H}^2} \\
& = 
\bar{V}_i^T (I + \mathbf{B}_{ii})^{-1/2}
\begin{bsmallmatrix}
\Phi_{i_1}U_{i_1}^{\mathcal{H}^2}  & & \\
& \ddots & \\
& & \Phi_{i_m}U_{i_m}^{\mathcal{H}^2} 
\end{bsmallmatrix}
R_i^{\mathcal{H}^2},
\end{align*}
which involves computations of $\{\Phi_{i_p}U_{i_p}^{\mathcal{H}^2}\}$ at level $(k-1)$. 

\paragraph{Case 2}
If we need the HSS coefficient matrix $B_{ij}$ and the 
corresponding $A_{ij}$ is a Type-1 block not at the leaf level, then 
$B_{ij}$ must be computed by recursion, using children blocks
in the $\mathcal{H}^2$ representation that are either Type-1 or Type-3.
If they are Type-3, then the recursion stops (we have the case above), but if they are Type-1, then
the recursion continues unless the Type-1 block is stored in dense format (i.e., at 
the leaf level, such a $B_{ij}$ is directly computed).

\paragraph{Case 3}
The case that we need the HSS coefficient matrix $B_{ij}$ and the
corresponding $A_{ij}$ is a Type-2 block is impossible (see the comment on 
Type-2 blocks in case 1); such $B_{ij}$ are never needed.

Overall, at each level $k$ of the accelerated SPD HSS construction, we only need to construct the small subset of all $B_{ij}$ blocks that are associated with either Type-1 or Type-3 blocks $A_{ij}$. 
There are in total only $O(|\text{lvl}(k)|)$ such blocks at level $k$.
Finally, we only require $B_{ij}$ for each pair of sibling nodes for the final HSS representation.
The computation of all such $B_{ij}$ is $O(Nr^2)$. 

\subsubsection{Summary} 
The complete algorithm that exploits an $\mathcal{H}^2$ representation to efficiently construct an SPD HSS approximation is shown in \Cref{alg:spdhss_h2}.
Note that only $B_{ij}$ corresponding to Type-3 and Type-1 blocks need to be computed.
The computation proceeds level-by-level from the leaves toward the root in
order to satisfy the data dependencies implicit in $B_{i j}$, $\mathbf{B}_{ii}$, and $\Phi_i$.

In the algorithm, the major computation and storage come from those related to $Y^{(k)}$. 
\Cref{alg:spdhss_h2} thus has $O(rN\log N)$ computation and peak storage cost.
The constructed SPD HSS approximation has $O(Nr)$ storage cost.

\cref{alg:spdhss_h2} can be extended to construct an SPD HSS approximation with a given approximation error threshold by adaptively adding more vectors to $\Omega$ and by compressing $\Lambda_{ii^c}\Omega_{i^c}$ with this error threshold. 
In this case, the approximation ranks for $\Lambda_{ii^c}\Omega_{i^c}$ could increase with the overall matrix size for many problems, leading to more expensive computation cost. 

\begin{algorithm}
\caption{Accelerated SPD HSS construction with quasilinear computation}
\label{alg:spdhss_h2}

\begin{algorithmic}
\REQUIRE HSS rank $r$, oversampling parameter $p$, an SPD $\mathcal{H}^2$ representation of  $A$
\ENSURE an SPD HSS approximation of $A$ with $\{A_{ii}\}, \{B_{ij}\}, \{U_i\}, \{R_i\}$

\STATE compute $Y^{(k)} = (A - \text{diag}(\{A_{ii}\}_{i\in\text{lvl}(k)}))\Omega, \quad k = 1,2,\ldots, L-1$

\textbf{At the leaf level}

\STATE \quad compute the Cholesky decomposition $A_{i i} = S_i S_i^T$, $\forall i \in \text{lvl}(1)$
		
\STATE \quad compute $\Lambda_{ii^c}\Omega_{i^c} = S_i^{-1}Y_i^{(1)}$, $\forall i \in \text{lvl}(1)$

\STATE \quad compute $V_i$ via the pivoted QR decomposition of $\Lambda_{ii^c}\Omega_{i^c}$ by \cref{alg:matvei_random}

\STATE \quad compute $\Phi_i U_i^{\mathcal{H}^2} =  V_i^TS_i^{-1} U_i^{\mathcal{H}^2}$, $\forall i \in \text{lvl}(1)$

\STATE \quad compute $B_{i j} = V_i^T S_i^{-1} A_{ij} S_j^{-T} V_j^T$ for all Type-1 $i,j\in \text{lvl}(1)$

\STATE \quad compute $B_{i j} = (\Phi_i U_i^{\mathcal{H}^2}) B_{ij}^{\mathcal{H}^2} (\Phi_j U_j^{\mathcal{H}^2})^T$ for all Type-3 $i,j\in \text{lvl}(1)$

\STATE \quad set $U_i = S_i V_i$, $\forall i \in \text{lvl}(1)$

\FOR {$k = 2,3,\ldots,L-1$}

\STATE compute $(I + \mathbf{B}_{ii})^{\pm 1/2}$ via eigen-decomposition of $\mathbf{B}_{ii}$, $\forall i \in \text{lvl}(k)$

\STATE compute $\Lambda_{ii^c}\Omega_{i^c} = (I+\mathbf{B}_{ii})^{-1/2} \left[\begin{smallmatrix}
\Phi_{i_1} & & \\
& \ddots & \\
& & \Phi_{i_m}
\end{smallmatrix}\right]Y_i^{(k)}$ via \Cref{alg:spdhss_matvec}, $\forall i \in \text{lvl}(k)$

\STATE compute $\bar{V}_i$ via the pivoted QR decomposition of $\Lambda_{ii^c}\Omega_{i^c}$ by \cref{alg:matvei_random}

\STATE compute $\Phi_i U_i^{\mathcal{H}^2} =  \bar{V}_i^T (I + \mathbf{B}_{ii})^{-1/2}
\begin{bsmallmatrix}
\Phi_{i_1}U_{i_1}^{\mathcal{H}^2} & & \\
& \ddots & \\
& & \Phi_{i_m}U_{i_m}^{\mathcal{H}^2}
\end{bsmallmatrix}
R_i^{\mathcal{H}^2}$, $\forall i \in \text{lvl}(k)$

\STATE compute $B_{i j} = \bar{V}_i^T (I + \mathbf{B}_{ii})^{-1/2} \mathbf{B}_{ij} (I + \mathbf{B}_{jj})^{-1/2}\bar{V}_j$,  for all Type-1 $i,j\in \text{lvl}(k)$

\STATE compute $B_{i j} = (\Phi_i U_i^{\mathcal{H}^2}) B_{ij}^{\mathcal{H}^2} (\Phi_j U_j^{\mathcal{H}^2})^T$ for all Type-3 $i,j\in \text{lvl}(k)$

\STATE set $R_i = (I+\mathbf{B}_{ii})^{1/2} \bar{V}_i$, $\forall i \in \text{lvl}(k)$
\ENDFOR
\end{algorithmic}
\end{algorithm}

%% file: doc/numerical_test.tex
\section{Numerical experiments}
\label{sec:numerical}

The SPD HSS approximation constructed by \cref{alg:spdhss_h2} will be denoted as ``SPDHSS.''
In comparison, the regular HSS representation that does not
use scaling-and-compression nor consider positive definiteness will
be referred to as ``regular HSS'' or, simply, ``HSS'' in the tables and figures below.

SPDHSS is tested using SPD kernel matrices.
A kernel matrix $K(X,X) = \left(K(x_i,x_j)\right)_{x_i,x_j\in X}$ 
is defined by a kernel function $K(x,y)$ and a set of points $X$. 
Kernel matrices appear in many applications, such as Gaussian processes 
and Brownian dynamics, and usually can be effectively represented in 
$\mathcal{H}^2$ form when defined in low-dimensional spaces, e.g., two-dimensional (2D) and three-dimensional (3D) spaces.
We consider four kernel functions:
\begin{itemize}
\item Mat\'{e}rn-$3/2$ kernel, $K(x,y) = \left(1+\sqrt{3}\,l\,|x-y|\right) \exp\left(-\sqrt{3}\,l\,|x-y|\right)$.
\item Gaussian kernel, $K(x,y) = \exp\left(-l \, |x-y|^2\right)$.
\item Inverse multiquadric (IMQ)  kernel, $K(x,y) = 1/\sqrt{1+ l \, |x-y|^2}$.
\item Rotne--Prager--Yamakawa (RPY) kernel \cite{rotne-prager,yamakawa},
\[ 
    K(x, y) = 
    \left\{ 
    \begin{array}{ll}
     \dfrac{1}{a}I_3 & \text{if } |r| = 0\\
     \dfrac{3}{4|r|} \left( I_3 + \dfrac{r r^T}{| r |^2}
     \right) + \dfrac{3a^2}{2| r |^3} \left( \dfrac{1}{3} I_3 - \dfrac{r r^T}{| r
     |^2} \right) & \text{if } |r| \geqslant 2a \\
      \dfrac{1}{a}\left( 1 - \dfrac{9}{32} \dfrac{| r |}{a}
     \right) I_3 + \dfrac{3}{32} \dfrac{| r |}{a} \dfrac{r r^T}{| r |^2}
     & \text{if } |r| < 2a
   \end{array} 
   \right.
\]
with $r = x-y$. The kernel is a $3\times 3$ tensor and is
defined for points in 3D.
\end{itemize}

The first three kernels are commonly used in statistical models with spatial data, such as Gaussian processes for geoscience problems \cite{heaton2019case},
as well as many other numerical methods that rely on radial
basis functions, such as in the numerical solution of partial differential
equations \cite{fornberg2015solving}.
In these kernels, $l$ is a length-scale parameter that is optimized to fit the data.
The RPY kernel describes the hydrodynamic interactions between spherical particles in a viscous fluid.  In this kernel, $a$ is the particle radius. 
In practice, the parameters $l$ and $a$ in these kernel functions and the distribution of the data points affect the conditioning of the resulting kernel matrices. 

For all tests in this section, we consider two types of 3D point sets for $X$: uniform random distributions of $N$ points on a 
sphere of radius $\sqrt{N/(4\pi)}$ in 3D (\textit{sphere} point set), and uniform random 
distributions of $N$ points in a ball of radius $\sqrt[3]{3N/(4\pi)}$ in 3D (\textit{ball} point set).
The radii of the sphere and the ball are selected to make the point density on the sphere and in the ball
remain constant with different $N$. 

We use the H2Pack library \cite{huang_h2pack_2020} for general computations 
related to $\mathcal{H}^2$ and HSS representations of kernel matrices. 
H2Pack can efficiently construct an $\mathcal{H}^2$ representation of a kernel matrix with 
linear-scaling computation by using a hybrid compression technique called the proxy point method \cite{xing2020interpolative}.
H2Pack also provides efficient \textit{regular} HSS construction for kernel matrices using the proxy point method.
This regular HSS construction \cite{xing2018efficient} does not use the scaling-and-compression technique and instead exploits analytic information of a kernel function to reduce the construction cost. 
It resembles recursive skeletonization \cite{ho_fast_2012} but 
works for general kernel matrices and requires an additional ULV decomposition for matrix
inversion. 
Note that the ULV decomposition has relatively cheap computation cost compared to the corresponding HSS construction. 
All timings of regular HSS and SPDHSS construction reported below include that for ULV decomposition.

Given a kernel matrix $K(X, X)$ with $X$ in $d$-dimensional space, a $2^d$-ary partition tree 
is constructed by recursively partitioning a box enclosing all the points $X$ (by bisecting
each dimension) until each finest box 
has less than 400 points.
This partition tree is used to construct the $\mathcal{H}^2$ representation of the kernel matrix.
The regular HSS and SPDHSS representations use the same partition tree.
The preconditioned conjugate gradient (PCG) method is used to solve kernel matrix systems.
The systems have random right-hand side vectors with entries chosen from the uniform distribution on $[-0.5, 0.5]$.  The PCG relative residual norm stopping threshold is $10^{-4}$.

The test calculations are carried out on a dual Intel Xeon Gold 6226 CPU computer with a total of 24 cores 
and 180 GB memory. One hyperthread per core is used. 
All codes are implemented in C and parallelized using OpenMP.

\subsection{Computational efficiency}
\begin{figure}[htbp]
\centering
\includegraphics[width=0.8\textwidth]{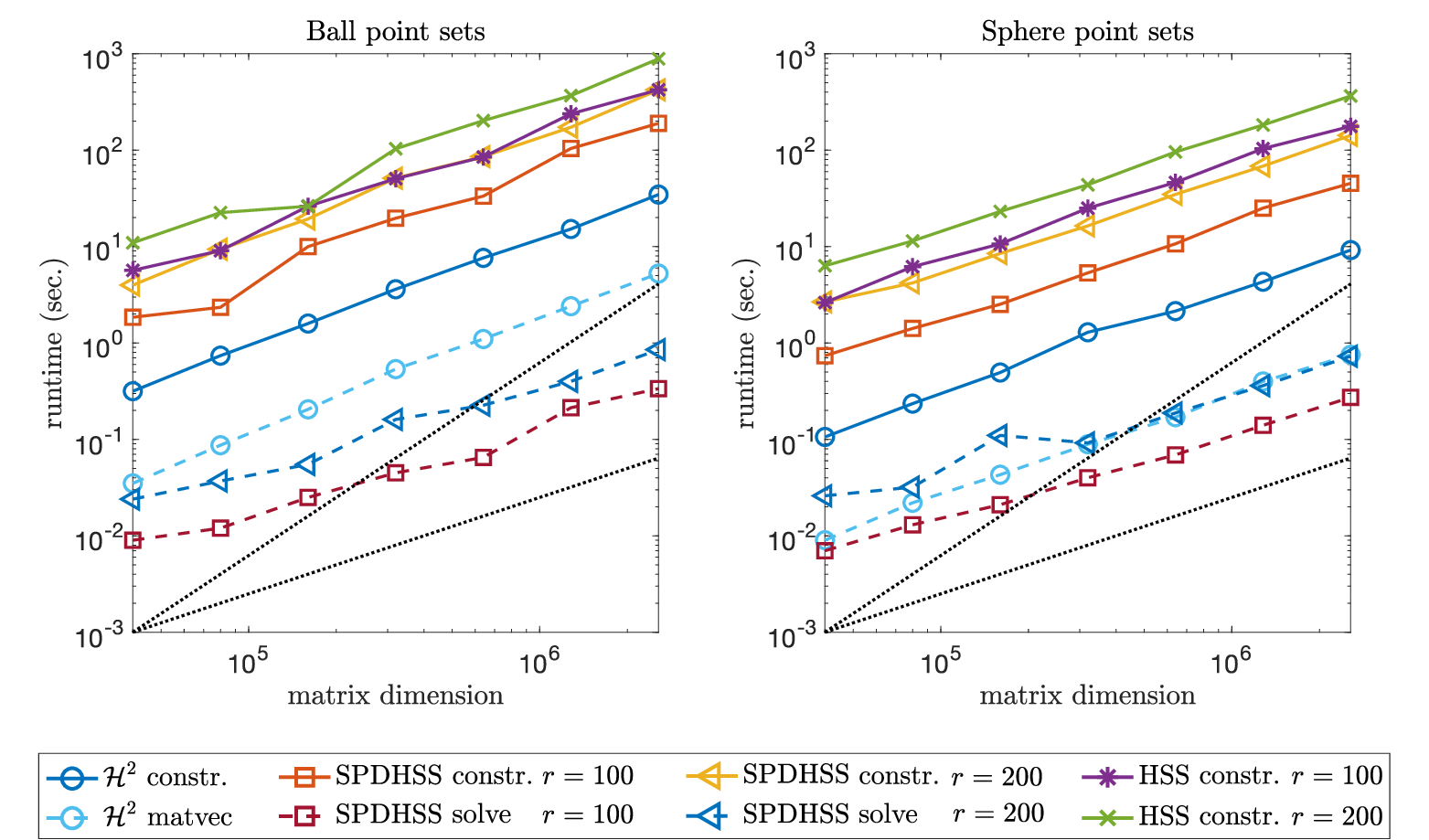}
\caption{Timings for constructing (solid lines) and applying (dashed lines) $\mathcal{H}^2$ representations and HSS/SPDHSS approximations.
Results are for Mat\'{e}rn kernel matrices with ball and sphere point sets. 
Linear-scaling of the regular HSS construction is due to the fixed approximation rank and the use of the proxy point method.
Dotted reference lines show linear and quadratic scaling. 
}\label{fig:test1_runtime}
\end{figure}

Consider the Mat\'{e}rn kernel with parameter $l=0.1$.
Ball and sphere point sets are generated with the number of points $N$ ranging from $4\times 10^4$
to $2.56\times 10^6$.
\cref{fig:test1_runtime} plots the timings for constructing $\mathcal{H}^2$, SPDHSS, and regular
HSS representations, as well as timings for $\mathcal{H}^2$ matrix-vector multiplication
and the SPDHSS/HSS solve operation (with a ULV decomposition).
The $\mathcal{H}^2$ representations are constructed with relative error threshold $10^{-8}$ here
and in the results that follow.  The SPDHSS/HSS approximations use fixed $r=100$ and $r=200$.

Approximate linear scaling is observed in all cases.
In terms of absolute cost, note that
SPDHSS construction formally requires $(r+p)(L-1)$ matrix-vector multiplications
using the $\mathcal{H}^2$ representation to compute $Y^{(k)}$ in \cref{eqn:Yk},
not to mention other operations.  
For our range of $N$, the number of levels $L$ ranges from 4 to 7.
The oversampling parameter $p$ is set to $10$ for all tests.
Despite requiring $(r+p)(L-1)$ matrix-vector multiplications,
\cref{fig:test1_runtime} shows that the SPDHSS construction time
can be faster than $r$ times the cost of a single $\mathcal{H}^2$ matrix-vector multiplication.
This is due to the use of level 3 BLAS operations when performing these
multiplications, and points to the computational efficiency
of blocked matrix multiplication that can be used in randomized algorithms.

\Cref{fig:test1_runtime} also shows that 
SPDHSS construction is faster than regular HSS construction with the same rank $r$.
This is due to efficient use of the $\mathcal{H}^2$ representation for SPDHSS construction
as presented in this paper.  Note that the $\mathcal{H}^2$ construction cost is relatively
very small.  We also note that the cost of the SPDHSS solves is comparable to or smaller
than the cost of a $\mathcal{H}^2$ matrix-vector multiplication in these examples.

\begin{figure}[h]
\centering
\includegraphics[width=0.8\textwidth]{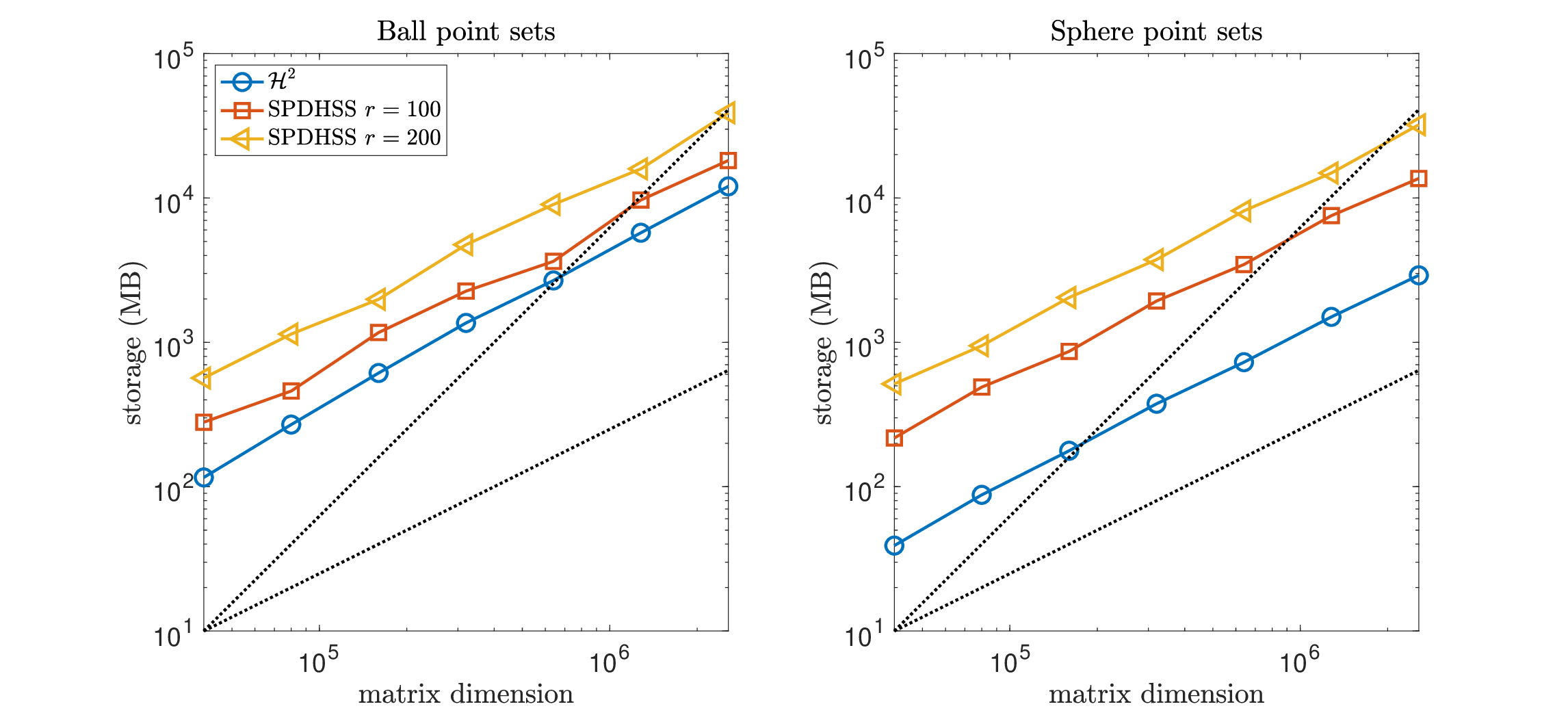}
\caption{Storage cost of $\mathcal{H}^2$ representations and SPDHSS approximations of Mat\'{e}rn kernel matrices with ball and sphere point sets. 
Dotted reference lines show linear and quadratic scaling. }\label{fig:test1_storage}
\end{figure}

\Cref{fig:test1_storage} plots the storage costs of the $\mathcal{H}^2$ representation
and SPDHSS approximation (after ULV decomposition).
As can be shown analytically \cite{xia_fast_2010}, an HSS approximation using a fixed rank has linearly scaling storage cost. 
An SPDHSS approximation has the same storage cost as the corresponding regular HSS approximation 
for the same approximation rank $r$.

\begin{figure}[h]
	\centering
	\includegraphics[width=0.8\textwidth]{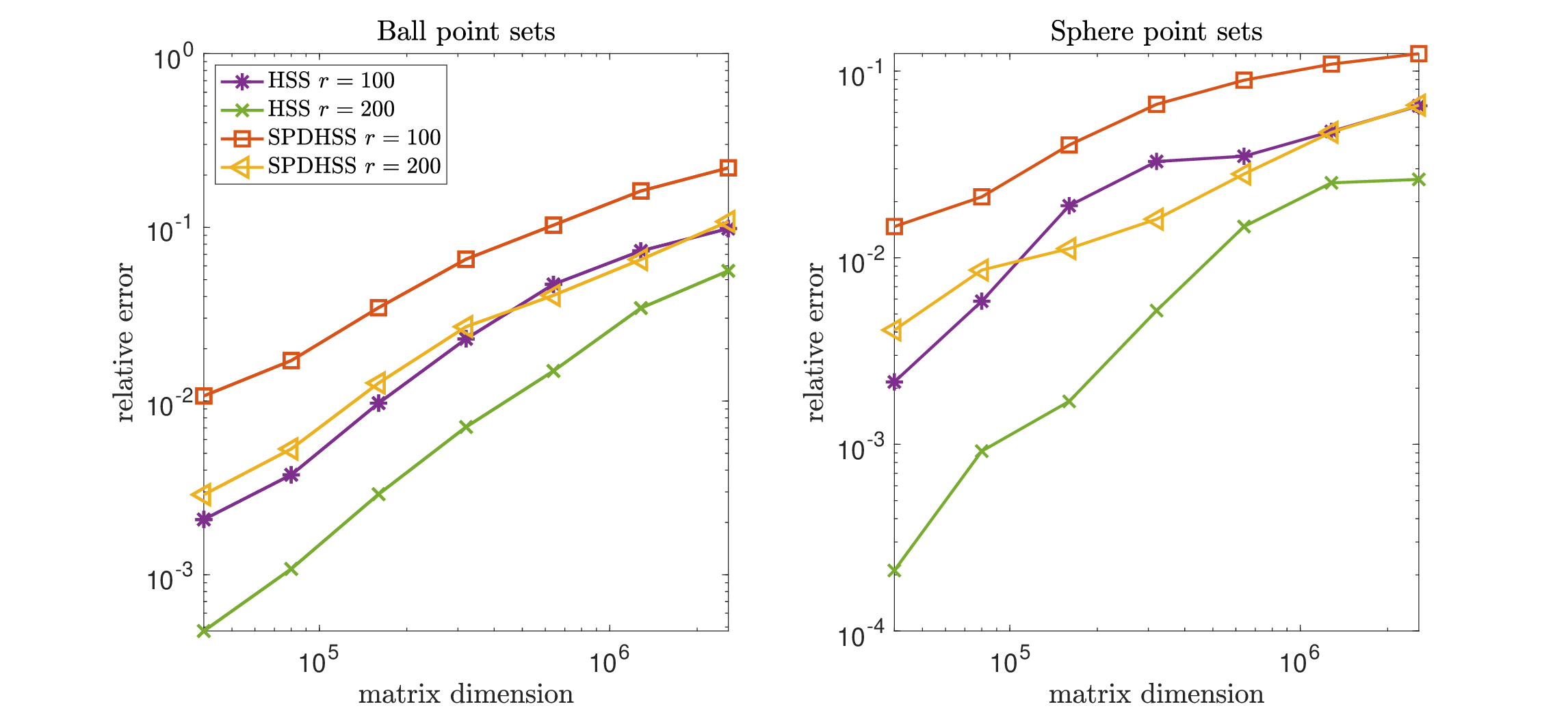}
	\caption{Average relative matrix-vector multiplication errors of regular HSS and SPDHSS
	approximations.
	Each data point is the average relative error (in 2-norm) of  matrix-vector multiplications by 10 Gaussian random vectors. 
    The reference results for matrix-vector multiplication are computed using the $\mathcal{H}^2$ representation.
	}\label{fig:test1_error}
\end{figure}

To estimate the accuracy of the regular HSS and SPDHSS approximations, we measure the accuracy of sample
matrix-vector multiplications with these approximations, where we assume that
the matrix-vector multiplication with the $\mathcal{H}^2$ representation is the exact value.
\Cref{fig:test1_error} plots the average relative error using a sample of 10
matrix-vector multiplications by Gaussian random vectors.
As expected in 3D problems, with a fixed rank $r$, the relative errors of both the SPDHSS and regular HSS approximations increase with the problem size.
SPDHSS has slightly larger approximation errors than regular HSS. 
However, the regular HSS approximations in all these examples are not SPD.

\subsection{Preconditioning performance}
We test the SPDHSS approximation as a preconditioner and compare it with the following preconditioners.
\begin{itemize}
\item The block Jacobi preconditioner (BJ) is a block diagonal matrix consisting of the diagonal blocks associated with the leaf nodes in the partition tree.
\item The factorized sparse approximate inverse preconditioner (FSAI) is $G^T G$, where $G$ is a sparse approximation to the inverted Cholesky factor of an SPD matrix \cite{kolotilina1993factorized}.
The nonzero pattern used for row $i$ of $G$ has $k$ nonzero entries and corresponds to the $k$ nearest neighbors of point $i$ for scalar kernels, or to the $k/3$ nearest neighbors of point $i$ for the RPY kernel.
Constructing $G$ requires only selected entries of the SPD matrix, which depends on the chosen
sparsity pattern, and thus FSAI can be efficient for dense kernel matrices.
\item The regular HSS approximation with a fixed rank $r$, but only when the approximation happens to be SPD.
\end{itemize}

\subsubsection{Kernel functions with varying parameters}
We consider the Mat\'{e}rn, Gaussian, and IMQ kernel functions with varying parameter $l$. 
In Gaussian process estimation, $l$ changes in each optimization step and each $l$ 
corresponds to a system to solve involving the kernel matrix, denoted here as $K_l(X, X)$.
For all three kernel functions, when $l$ is close to zero, $K_l(X, X)$ is close to 
low-rank; 
when $l$ is sufficiently large, $K_l(X,X)$ is close to sparse. 
In practice, a diagonal shift is added, i.e., $\sigma I + K_l(X,X)$, to 
account for noise in the Gaussian process model.
Numerically, this diagonal shift is also necessary to keep the linear system
from being extremely ill-conditioned when $l$ is small. 
We set $\sigma = 10^{-2}$ in the following tests.

\Cref{tab:iternum} lists the number of PCG iterations for solves involving 
the three types of kernel matrices
with various parameters $l$, generated by a ball point set of size $ 3.2\times 10^5$. 
For all the different types of kernels, FSAI performs very well when $l$ is large, corresponding to
kernel matrices that are close to sparse.
For smaller values of $l$, the performance of FSAI deteriorates.

In comparison, SPDHSS has more consistent preconditioning performance for this wide range of parameters, although it takes more iterations than FSAI for large $l$ in most cases. 
This consistency is an advantage of SPDHSS over FSAI,
since the parameter $l$ changes during optimization,
and it could be difficult to quantitatively decide when to use FSAI,
particularly for more complicated kernel functions.

We also observe that SPDHSS has better performance than BJ in all the tests. 
The SPDHSS preconditioner can be viewed as the combination of a BJ preconditioner (the diagonal blocks) with some off-diagonal approximations. 
Finally, the computed regular HSS preconditioner in most cases is not SPD.

\begin{table}[hp]
	\centering
\footnotesize
	\caption{Number of PCG iterations for systems with diagonal-shifted kernel matrices $\sigma I+ K_l(X, X)$ with different kernels and parameters $l$.
All tests use the same $N = 3.2\times10^5$ points in a ball.
The notation ``$-$'' means that PCG fails to converge within $3000$ iterations; ``/'' means that a regular HSS preconditioner is not SPD. 
	}
	\label{tab:iternum}
	\subfloat[Mat\'{e}rn kernel]{
		     	\setlength{\tabcolsep}{3.6pt}
				\begin{tabular}{l| rrrrrrrrrr}
					\toprule
parameter $l$ & 0.0010 & 0.0025 & 0.005 & 0.010 & 0.025 & 0.05 & 0.10 & 0.25 & 0.5 & 1.0 \\
\midrule
Unpreconditioned & 41 & 119 & 297 & 687 & 1896 & - & - & 1684 & 634 & 210 \\
BJ & 459 & 1202 & 2504 & - & - & - & 2765 & 707 & 172 & 82 \\
FSAI $k=200$ & 2659 & 2623 & 2045 & 1518 & 960 & 569 & 266 & 63 & 18 & 6 \\
FSAI $k=400$ & 1734 & 1531 & 1111 & 831 & 511 & 266 & 108 & 28 & 10 & 4 \\
SPDHSS $r=100$ & 2 & 3 & 5 & 14 & 62 & 159 & 287 & 245 & 117 & 56 \\
SPDHSS $r=200$ & 1 & 2 & 4 & 5 & 20 & 58 & 131 & 171 & 92 & 50 \\
HSS $r=100$ & 2 & 3 & 6 & / & / & / & / & / & / & / \\
HSS $r=200$ & 2 & 2 & 3 & 6 & / & / & / & / & / & / \\
				\bottomrule
			\end{tabular}
		}
	
		\subfloat[Gaussian kernel]{
		\setlength{\tabcolsep}{5pt}
		\begin{tabular}{l| rrrrrrrrr}
			\toprule
parameter $l$ & 0.0001 & 0.0005 & 0.001 & 0.005 & 0.01 & 0.05 & 0.1 & 0.5 & 1.0 \\
\midrule
Unpreconditioned & 103 & 394 & 838 & 2498 & 2396 & 1098 & 695 & 232 & 150 \\
BJ & 1851 & - & - & - & - & 989 & 487 & 147 & 85 \\
FSAI $k=200$ & - & - & - & 2718 & 1371 & 253 & 106 & 16 & 7 \\
FSAI $k=400$ & - & - & - & 1498 & 687 & 93 & 43 & 9 & 4 \\
SPDHSS $r=100$ & 1 & 3 & 21 & 577 & 855 & 563 & 378 & 118 & 65 \\
SPDHSS $r=200$ & 1 & 2 & 6 & 116 & 284 & 415 & 314 & 106 & 59 \\
HSS $r=100$ & 2 & 3 & / & / & / & / & / & / & / \\
HSS $r=200$ & 2 & 2 & 2 & / & / & / & / & / & / \\
			\bottomrule
		\end{tabular}
	}

	\subfloat[IMQ kernel]{
	\setlength{\tabcolsep}{3.6pt}
	\begin{tabular}{l| rrrrrrrrrrr}
		\toprule
parameter $l$ & 0.001 & 0.005 & 0.01 & 0.05 & 0.1 & 0.5 & 1 & 5 & 10 & 50 & 100 \\
\midrule
Unpreconditioned & 1239 & 2656 & - & - & 2812 & 1958 & 1576 & 915 & 724 & 394 & 284 \\
BJ & - & - & - & 2300 & 1605 & 529 & 322 & 195 & 151 & 97 & 73 \\
FSAI $k=200$ & 2839 & 1092 & 619 & 201 & 121 & 45 & 32 & 24 & 24 & 23 & 23 \\
FSAI $k=400$ & 1598 & 535 & 266 & 73 & 51 & 24 & 21 & 18 & 18 & 18 & 17 \\
SPDHSS $r=100$ & 53 & 249 & 328 & 355 & 291 & 112 & 74 & 38 & 28 & 13 & 10 \\
SPDHSS $r=200$ & 12 & 71 & 129 & 212 & 195 & 84 & 60 & 33 & 23 & 10 & 8 \\
HSS $r=100$ & / & / & / & / & / & / & / & / & / & / & / \\
HSS $r=200$ & / & / & / & / & / & / & / & / & / & / & / \\
		\bottomrule
	\end{tabular}
   }
\end{table}

\begin{table}[htbp]
\footnotesize
	\caption{
Timings (in sec.) for constructing and applying (solves) the preconditioner and
storage (in GB) for the BJ, FSAI, and SPDHSS preconditioners.
The table also includes timings for constructing and applying (matrix-vector multiplication)
the $\mathcal{H}^2$ representation.  
The matrices are from those in \Cref{tab:iternum} defined by the Mat\'{e}rn with $l=0.025$, Gaussian with $l=0.01$, and IMQ with $l=1$.
The $\mathcal{H}^2$ representations of the three test matrices require $0.7$, $1.8$, and $1.2$ GB of storage, respectively. 
	}\label{tab:info}
	\centering
	\begin{tabular}{l|c|rl|rl|rl}
		\toprule
		          &     & \multicolumn{2}{c|}{Mat\'{e}rn} & \multicolumn{2}{c|}{Gaussian} & \multicolumn{2}{c}{IMQ}\\
		 & storage & constr.  & apply  & constr.  & apply   & constr.  & apply  \\
		\midrule
		$\mathcal{H}^2$ representation  & 0.7/1.8/1.2 &  1.7   &  0.30     & 10.8  &  0.36   &  4.6 &  0.92   \\  
		BJ                              & 0.3         &  0.3   &  0.016    &  2.9  &  0.0055 &  0.1 &  0.0091 \\  
		FSAI $k = 200$                  & 0.7         &  8.6   &  0.0070   &  6.7  &  0.0052 &  7.8 &  0.0083 \\  
		FSAI $k = 400$                  & 1.4         & 13.1   &  0.013    & 12.3  &  0.0081 & 15.3 &  0.023  \\  
		SPDHSS $r=100$                  & 2.2         & 21.1   &  0.046    & 22.5  &  0.041  & 23.6 &  0.043  \\  
		SPDHSS $r=200$                  & 4.6         & 38.3   &  0.15     & 44.0  &  0.11   & 43.9 &  0.11  \\  
		\bottomrule
	\end{tabular}
\end{table}

\Cref{tab:info} shows the time required to construct and apply (solve with) the various preconditioners.
The table also shows the time required to construct and apply (multiply by) the $\mathcal{H}^2$ representation.
The storage requirements for the $\mathcal{H}^2$ representation and for the preconditioners are also shown.
The construction cost of an SPDHSS approximation depends on the efficiency of the corresponding $\mathcal{H}^2$ representation, and thus varies for different kernel functions.
The application of the SPDHSS preconditioners, although more expensive than for
FSAI preconditioners, is relatively fast in comparison to corresponding 
$\mathcal{H}^2$ matrix-vector multiplications.

\subsubsection{Kernel matrices with varying sizes}
We now consider the iterative solution of the Mat\'{e}rn and RPY kernel matrix systems for systems of different sizes.
For the RPY kernel, particle radii $a = 0.29$ and $a = 0.42$ are selected such that 
each ball point set has corresponding volume fraction of particles around 0.1 and 0.3, respectively.
These two volume fractions are representative for macromolecular simulations of conditions
within biological cells \cite{chow2016pnas}.
For the Mat\'{e}rn kernel, $l= 0.25$ and $l=0.01$ are tested based on the previous results in \cref{tab:iternum},
where FSAI performs better than SPDHSS for $l=0.25$ and vice versa for $l=0.01$. 
No diagonal shift is added to RPY kernel matrices while a shift of $\sigma = 10^{-2}$ is added to Mat\'{e}rn kernel matrices as before.

\Cref{tab:iternum1} shows PCG convergence for systems using the two kernel functions with different point sets. 
As expected, iteration numbers increase with matrix sizes for the FSAI and SPDHSS preconditioners since a fixed approximation rank $r$ and sparsity parameter $k$ are used. 
Related to this is the increasing relative approximation error in the SPDHSS approximation with increasing matrix size when $r$ is fixed, as observed earlier in \cref{fig:test1_error}.
As to be shown next, it is possible to obtain scalable preconditioning performance by constructing SPDHSS preconditioners with a fixed relative error threshold but at the sacrifice of asymptotically more expensive cost in SPDHSS construction and solve.

\begin{table*}[htbp]
\footnotesize
\caption{
Number of PCG iterations for kernel matrices defined by different point sets. 
}\label{tab:iternum1}
\centering
\setlength{\tabcolsep}{3.6pt}
\begin{tabular}{l|ccccc|ccccc}
\toprule
\multicolumn{1}{c|}{\multirow{2}{*}{$N$ ($\times 10^4$)}}	& \multicolumn{5}{c|}{ball point sets} & \multicolumn{5}{c}{sphere point sets} \\
 & 4 & 8 & 16 & 32 & 64 & 4 & 8 & 16 & 32 & 64 \\
 
\midrule
\textit{Mat\'{e}rn $l=0.01$} & & & & & & & & & & \\
Unpreconditioned        &  134   & 217     & 397      & 689    & 1235  & 370   & 713    & 1252 & 2195   & - \\  
BJ                                &  1209 &  2249 & -             & -          &  -          & 1472 & 2642 &  -         & -             & - \\
FSAI $k = 200$      &  762   &  1215 &  1205   & 1503 &  -          & 504   & 651     & 877  & 1340    & 1791 \\
FSAI $k = 400$      &  347   &  693 &   656       &  849   &  1474 & 247   & 272    &  380  & 470      & 634 \\
SPDHSS $r=100$  &  3        &  4       &   7            &  14      &   29     & 4        & 4          & 7        & 17         & 34 \\
SPDHSS $r=200$  &  2        &  3        &  3           &   5         &   10     & 3        & 3         &  3       & 4            & 8 \\
 

 \midrule  
\textit{Mat\'{e}rn $l=0.25$} & & & & & & & & & & \\
Unpreconditioned        &   1123 &  1358  &   1542    &  1681  & 1790  & 549  & 559  & 566  & 569 & 570 \\
BJ                                &   574    &  568     &   763      &   704    & 650     & 169  & 197  & 184  &199  & 197 \\
FSAI $k = 200$      &   56       &  83       &   53         &   63       & 114    & 16     & 23     & 21     & 24   &  28 \\
FSAI $k = 400$      &   26       &  36       &   22         &   28       & 46       & 7       & 9       &  8       & 10    & 12 \\
SPDHSS $r=100$  &   100     &  148    &   196      &   236     & 294    & 38    & 55     &  71    & 85    & 101 \\
SPDHSS $r=200$  &   45        &  76      &   122      &   172     & 216     & 11   & 21     &  36    & 51    & 71 \\

\midrule
\textit{RPY $a=0.29$} & & & & & & & & & & \\
Unpreconditioned        & 432    & 510 & 1055 & 1151 & 1653 & 549 & 707  & 1741 & 1448 & 1706 \\  
BJ                                & 95    & 142 & 181 &  238 &  282     & 51 & 104  & 100 & 139 & 147 \\
FSAI $k = 200$      & 89    & 113 & 137 &  174 &  222     & 31 & 35  & 42 & 50 & 60 \\
FSAI $k = 400$      & 84    & 98 & 122 &  162 &  212     & 24 & 28  & 34 & 37 & 48  \\
SPDHSS $r=100$  & 40    & 51 & 62 &  77 &  99     &  16   & 22  & 27 & 28 & 31\\
SPDHSS $r=200$  & 27    & 36 & 46 &  58 &  75     &  12   & 14  & 22 & 22 & 25\\
\midrule  
\textit{RPY $a=0.42$} & & & & & & & & & & \\
Unpreconditioned         & 632 & 762 & 1571 & 1723 & 2428   & 809 & 1063 & 2566 & 2132 & 2515\\  
BJ                                 & 150 &  218   &  237 &  328  &  436      & 66 & 140 & 128 & 178 & 191\\
FSAI $k= 200$        & 137 &  165   &  206 &  266  &  338      & 38 & 43 & 52 & 62 & 77\\
FSAI $k= 400$        & 125 &  151   &  176 &  244  &  307      & 28 & 32 & 37 & 43 & 56\\
SPDHSS $r=100$   & 64 &  79   &  101 &  121       &  157      & 23 & 30 & 35 & 39 & 45\\       
SPDHSS $r=200$   & 45 &  59   &  74 &  95            &  122      & 16 & 20 & 30 & 29 & 36\\                                  
\bottomrule
\end{tabular}
\end{table*}

\subsubsection{SPDHSS with a fixed relative error threshold}
To demonstrate the preconditioning performance and computational complexity of SPDHSS with a fixed relative error threshold, we consider the Mat\'{e}rn kernel with $l=0.25$ which was previously tested with fixed ranks (\cref{tab:iternum1}).
Applying SPDHSS with two relative error thresholds $\tau = 10^{-1}$ and $\tau = 10^{-2}$, \Cref{tab:test_reltol} shows PCG iteration counts and average matrix-vector multiplication errors of SPDHSS approximations for systems with different point sets.
The iteration counts are roughly constant with different-sized problems, suggesting scalable preconditioning performance.
The results in \Cref{tab:test_reltol} show that the error of an SPDHSS approximation is well controlled by the relative error threshold for our test problems.
However, we note that an error threshold is applied to the compression of scaled blocks in SPDHSS construction, which only indirectly controls the overall matrix approximation error.

\begin{table*}[htbp]
	\footnotesize
	\caption{
		Number of PCG iterations and average matrix-vector multiplication errors  for kernel matrices defined by the Mat\'{e}rn kernel with $l=0.25$ and different point sets. 
	\label{tab:test_reltol}	
}
	
	\centering
	\setlength{\tabcolsep}{3.6pt}
	\begin{tabular}{l|ccccc|ccccc}
		\toprule
		\multicolumn{1}{c|}{\multirow{2}{*}{$N$ ($\times 10^4$)}}	& \multicolumn{5}{c|}{ball point sets} & \multicolumn{5}{c}{sphere point sets} \\
		& 4 & 8 & 16 & 32 & 64 & 4 & 8 & 16 & 32 & 64 \\
		
		\midrule
		\textit{Iteration\ counts} & & & & & & & & & & \\
		SPDHSS $\tau=10^{-1}$  &  85   &  79   &  90   &  88  &   89   &  43   &  43   &  46  &  50  &  49 \\
		SPDHSS $\tau=10^{-2}$  &  10   &  10   &  10   &  13  &   10      & 7      & 9     &  12      &   10  & 10 \\
		
		\midrule
		\textit{Matvec errors} & & & & & & & & & & \\
		SPDHSS $\tau=10^{-1}$  &  0.05  &  0.05  &  0.06  &  0.06   &  0.06    & 0.06        & 0.06          & 0.06       & 0.06         & 0.06 \\
		SPDHSS $\tau=10^{-2}$  &  0.004  &  0.003  &  0.003  &  0.003   &   0.003   & 0.005     & 0.004         &  0.005       & 0.005   & 0.005\\
		\bottomrule
	\end{tabular}
\end{table*}

\begin{figure}[htbp]
	\centering
	\subfloat[Maximum ranks]{
		\includegraphics[width=0.45\textwidth]{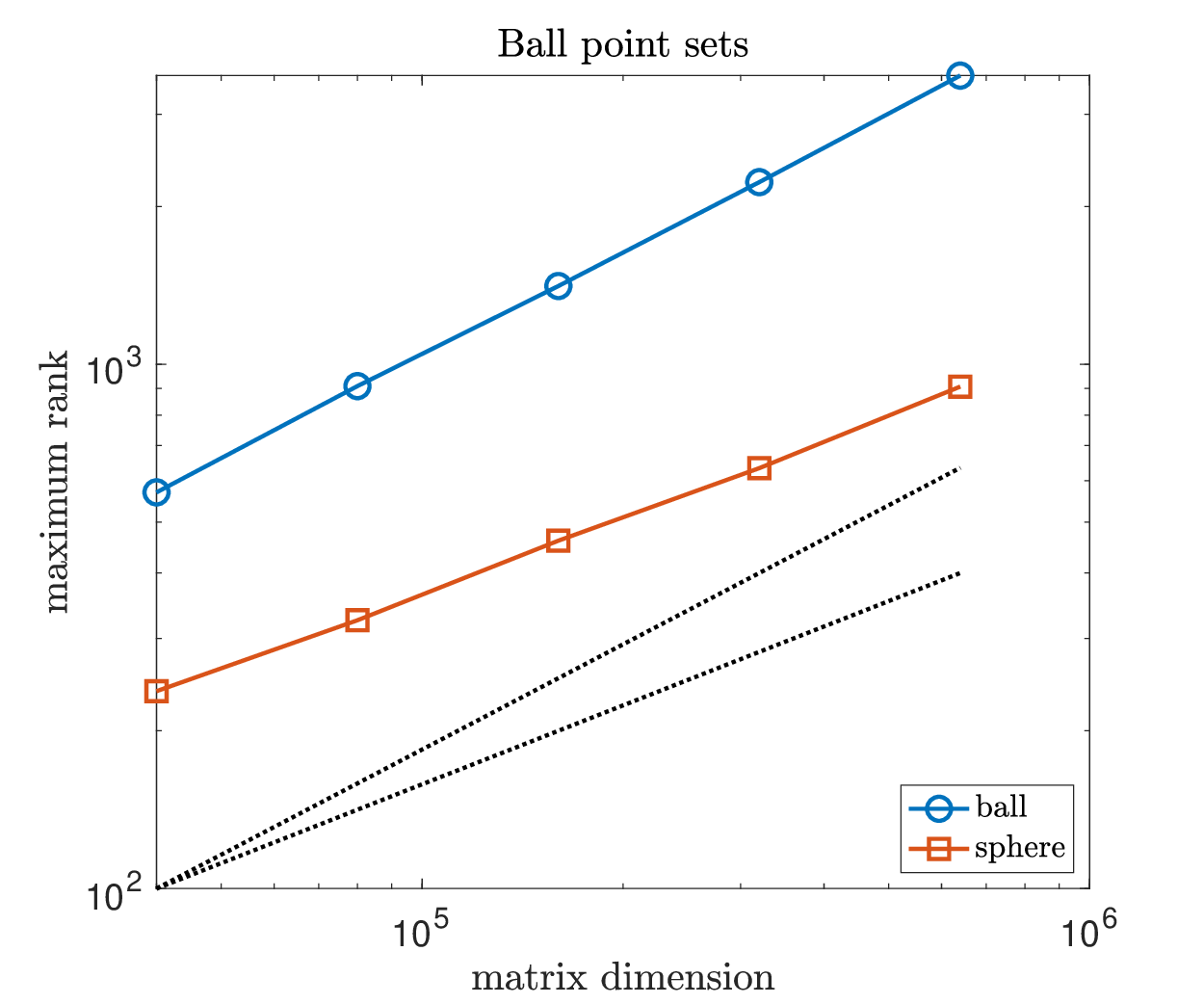}
	}
	\subfloat[Timings for construction and application]{
		\includegraphics[width=0.45\textwidth]{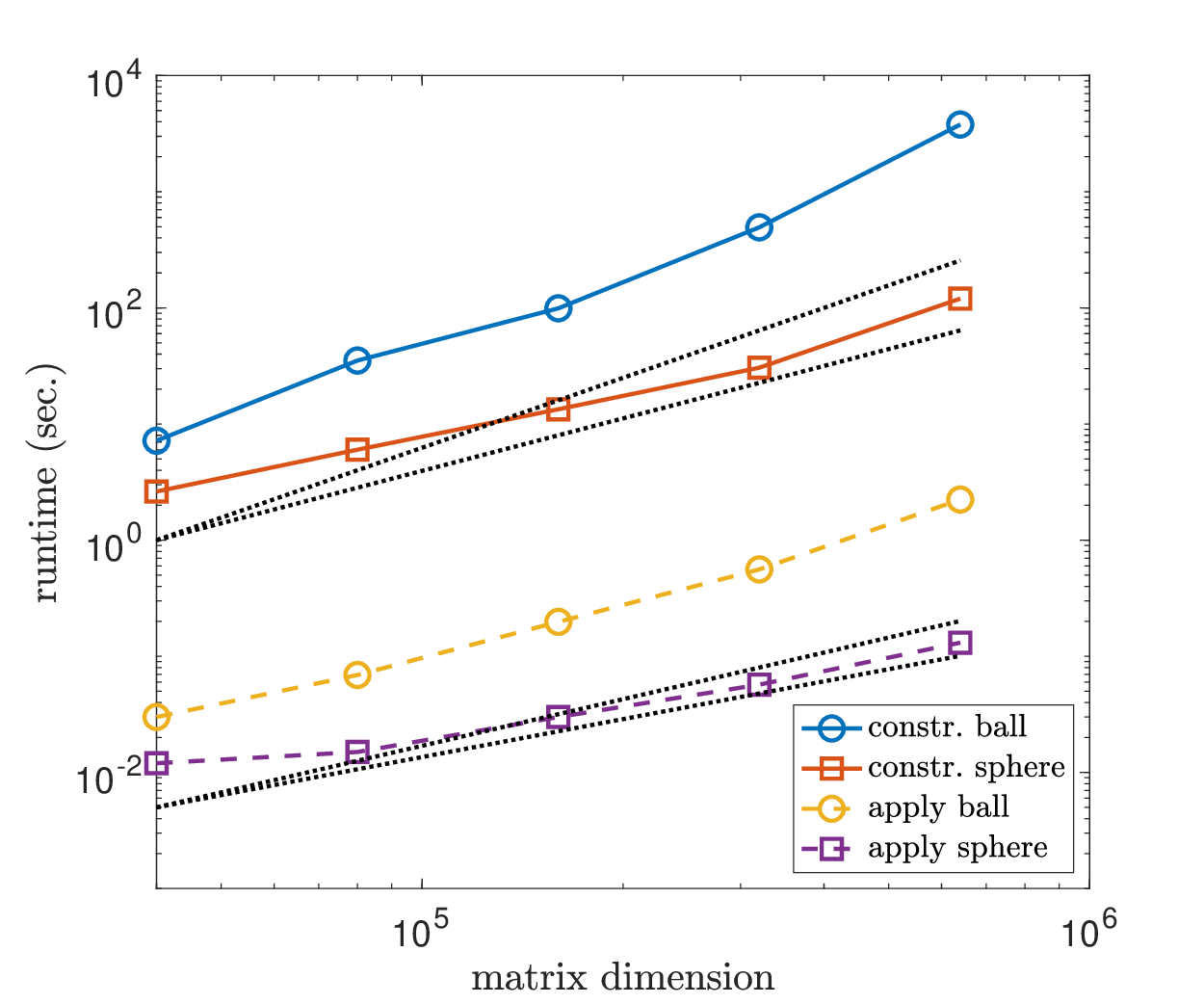}
	}
	\caption{
	Maximum approximation ranks and timings for construction and application (solves) of SPDHSS with fixed relative error threshold $\tau = 10^{-2}$ for ball and sphere point sets. 	
	In (a), the two dotted reference lines show $O(N^{1/2})$ and $O(N^{2/3})$ scaling. 
	In (b), the four dotted reference lines show $O(N\log(N))$, $O(N^{4/3})$, $O(N^{3/2})$, and $O(N^2)$ scaling. 
	}
	\label{fig:test_reltol}
\end{figure}

\Cref{fig:test_reltol} plots the maximum approximation ranks and SPDHSS construction and application costs with different point sets for $\tau = 10^{-2}$.
For these ball (sphere) point sets, it is well known (e.g., see \cite{ho_fast_2012}) that an HSS approximation with a fixed relative error, whether it is regular HSS or SPDHSS, can have at least $O(N^{2/3})$ ($O(N^{1/2}$)) maximum approximation ranks and have at least $O(N^2)$ ($O(N^{3/2})$) factorization and $O(N^{4/3})$ ($O(N\log N)$) solve cost.  
Our numerical results with SPDHSS are consistent with this theoretical analysis.

Overall, although the number of PCG iterations remains roughly constant when the problem size increases, the costs of precomputation of SPDHSS construction and the application of SPDHSS as a preconditioner both increase superlinearly. It is thus more practical to use a proper combination of a maximum rank threshold and a relative error threshold for the application of SPDHSS.

%% file: doc/conclusion.tex
\section{Conclusion}

Fast direct solvers and rank-structured preconditioners, such as those using
the HSS representation, impose a block structure on a matrix that provides
for fast solve operations, but the rigid block structure (arising from
so-called ``weak admissibility'') also results in large block ranks,
especially if an accurate representation is desired.  This leads to high
construction cost.

On the other hand, more general rank-structured matrix representations,
such as $\mathcal{H}^2$, have a flexible block structure (arising
from so-called ``strong admissibility'') that allows for an accurate
representation with smaller block ranks, and thus these representations
have relatively low construction cost.  However, the general structure
does not admit fast solve operations.

This paper, in a way, combines these two types of rank-structured matrix
representations.  The paper shows how to accelerate the construction
of an SPD HSS approximation to an SPD matrix by exploiting and only using an
$\mathcal{H}^2$ representation of the SPD matrix that is assumed to be
available, for example, in the context of a preconditioned iterative
solve.  The acceleration results from (i) using fast $\mathcal{H}^2$
matrix-vector multiplication to compute scaled basis matrices, $V_i$
and $\bar{V}_i$, needed in constructing the HSS representation, and
from (ii) using existing low-rank blocks in the $\mathcal{H}^2$
representation to reduce the number of coefficient matrices $B_{ij}$
that need to be computed in the HSS representation.

While we only tested SPDHSS as a preconditioner on kernel matrices, its
application to linear systems from the numerical solution of integral
equations is straightforward.  Further, although we only considered
dense SPD matrices, our proposed algorithms, \cref{alg:spdhss} and
\cref{alg:spdhss_h2}, can also be applied directly to sparse SPD matrices
that are ubiquitous in the numerical solution of partial differential
equations.  Both methods can still guarantee the positive definiteness
of the constructed preconditioners in the sparse case, but it is worthy
to study whether it is possible to exploit matrix sparsity directly to
accelerate SPD HSS construction.  It is also worthy to study
whether or not the FSAI and SPDHSS preconditioners can be beneficially
combined, i.e., augmenting a sparse preconditioner with a dense one.